\magnification=\magstep1   
\input amstex
\UseAMSsymbols
\input pictex 
\vsize=23truecm
\NoBlackBoxes
\parindent=18pt
  
   \font\rmk=cmr8


\def\Hom{\operatorname{Hom}}

\def\Ext{\operatorname{Ext}}

\def\rad{\operatorname{rad}}
\def\add{\operatorname{add}}
\def\Ker{\operatorname{Ker}}

\def\soc{\operatorname{soc}}
\def\Tr{\operatorname{Tr}}

\def\top{\operatorname{top}}
\def\col{\,\mathpunct{:}}
  \def\ss{\ssize }
\def\arr#1#2{\arrow <1.5mm> [0.25,0.75] from #1 to #2}
\def\arra#1#2{\arrow <2mm> [0.25,0.75] from #1 to #2}
\def\s{\hfill \square} 
\noindent
\vglue1truecm
\centerline{\bf  Gorenstein-projective and}
                     \smallskip
\centerline{\bf  semi-Gorenstein-projective modules. II}
                     \bigskip
\centerline{Claus Michael Ringel, Pu Zhang}
                \bigskip\medskip
 
\noindent {\narrower Abstract:  \rmk Let $\ss k$ be a field and 
$\ss q$ a non-zero element
of $\ss k$. In Part I, we have exhibited a
6-dimensional $\ss k$-algebra $\ss \Lambda = \Lambda(q)$ and we have shown that if 
$\ss q$ has infinite multiplicative order, then 
$\ss \Lambda$ has a 3-dimensional local
module which is semi-Gorenstein-projective, but not torsionless, thus not
Gorenstein-projective. This Part II is devoted to a detailed study of all the 
3-dimensional local $\ss \Lambda$-modules for this particular algebra $\ss\Lambda$.
If $q$ has infinite multiplicative order, we will encounter a whole family of 
3-dimensional local modules which are semi-Gorenstein-projective, but not torsionless.
	\bigskip
\noindent
Key words and phrases. Gorenstein-projective module, semi-Gorenstein-projective module, torsionless module, extensionless module, reflexive module, $\ss t$-torsionfree module, $\ss \mho$-quiver.
	\medskip
\noindent
2010 Mathematics Subject classification. Primary 16G10, 16G50. 
Secondary 16E05, 16E65, 20G42.
	\medskip
\noindent
Supported by NSFC 11431010
\par}

	\bigskip
{\bf 1. Introduction.}
	\medskip
{\bf (1.1)} We refer to our previous paper [RZ1] as Part I.
As in Part I, 
let $k$ be a field, and $q$ a non-zero element of $k$. We consider again the $k$-algebra
$\Lambda = \Lambda(q)$ generated by $x,y,z$ with relations
$$
 x^2,\ y^2,\ z^2,\ yz,\ xy+qyx,\ xz-zx,\ zy-zx.
$$
The algebra $\Lambda$ is a 6-dimensional local algebra 
with basis $1,x,y,z,yx,zx$. Its socle is $\soc\Lambda = \rad^2 \Lambda = 
\Lambda yx \oplus \Lambda zx.$ 
If not otherwise stated, all the modules considered 
will be left $\Lambda$-modules.

We follow the terminology used in Part I. In particular, we denote by $\mho M$
the cokernel of a minimal left $\add(\Lambda)$-approximation of $M$. 
In addition, we introduce the following definitions. 
We say that a module $M$ is {\it extensionless} if $\Ext^1(M,\Lambda) = 0$.
An indecomposable
semi-Gorenstein-projective module will be said to be {\it pivotal}
provided it is not torsionless. 
An indecomposable $\infty$-torsionfree
module will be said to be {\it pivotal} provided it is not extensionless.
Thus, a module $M$ is semi-Gorenstein-projective if and only if $\Omega^t M$
is extensionless for all $t\ge 0$; a torsionless module $M$ is
reflexive if and only if $\mho M$ is torsionless (see Part I (2.4));
a module $M$ is $\infty$-torsionfree if and only if $\mho^t M$ is reflexive
for all $t\ge 0$; and $M$ is Gorenstein-projective if and only if $M$ is both
semi-Gorenstein-projective and $\infty$-torsionfree.

	\medskip
{\bf (1.2)} We are interested in the semi-Gorenstein-projective and the 
$\infty$-torsionfree modules and will exhibit those which are
3-dimensional. We recall that a finite length
module is said to be {\it local} provided its top is simple. 
Thus, a local module is indecomposable; and if $R$ is a left artinian ring, then a
left $R$-module $M$ is local if and only if $M$ is a quotient of an indecomposable
projective module. A consequence of our study is the following assertion
	\medskip
{\bf Proposition.} {\it Let $M$ be a non-zero module of dimension
at most $3$. If $M$ is semi-Gorenstein-projective, then all the modules $\Omega^tM$
with $t\ge 0$ are $3$-dimensional and local. 
If $M$ is $\infty$-torsionfree, then all the modules $\mho^tM$
with $t\ge 0$ are $3$-dimensional and local.} In particular, if $M$ is
Gorenstein-projective, then all the modules $\Omega^tM$ and $\mho^t M$ 
with $t\ge 0$ are $3$-dimensional and local.
	\medskip
{\bf (1.3)} 
The text restricts the attention to the 3-dimensional local modules. 
The starting point of our investigation are two observations. The first one:
	\medskip
{\bf Proposition 1.} {\it A module of dimension at most $3$ is annihilated by $\rad^2\Lambda$,}
thus it is a module of Loewy length at most $2$.
	\medskip
The second observation is:
	\medskip
{\bf Proposition 2.} 
{\it An indecomposable $3$-dimensional torsionless module is local.}
	\medskip
The proof of Proposition 1 will be given in (2.6), the proof of Proposition 2 in (2.7).
	\medskip
{\bf (1.4) The  3-dimensional local modules.} 
We identify $(a,b,c)\in k^3\setminus\{0\}$ with $ax+by+cz$ 
and denote by $(a\col b\col c)$ the 1-dimensional subspace of $k^3$ generated by $(a,b,c)$.
The left ideal
$$
 U(a,b,c) = U(a\col b\col c) = \Lambda(a,b,c)+\soc \Lambda
$$
has dimension 3, and we obtain the left $\Lambda$-module
$$
 M(a,b,c) = M(a\col b\col c) = {}_\Lambda\Lambda/U(a,b,c).
$$
Clearly, {\it $M(a,b,c)$ is a $3$-dimensional local module and 
the modules 
$M(a,b,c)$, $M(a',b',c')$ are isomorphic if and only if  $(a\col b\col c) = (a'\col b'\col c').$}  
Let us add that the definition of $M(a,b,c)$ implies that  
$\Omega M(a,b,c) \simeq  U(a,b,c)$, this will be used throughout the text.

Conversely, {\it any $3$-dimensional local module is isomorphic to a module of the form
$M(a,b,c)$.}
In order to see this, one should look at the factor algebra 
$\overline\Lambda$ of $\Lambda$ modulo $\soc \Lambda = \rad^2\Lambda$,
thus $\overline \Lambda$ is the $k$-algebra generated by $x,y,z$ with relations all monomials
of length $2$. 
The $\Lambda$-modules of Loewy length at most 2 are just the modules annihilated 
by all monomials of length $2$, thus  the $\overline \Lambda$-modules. 
It is clear that the modules $M(a,b,c) = \overline\Lambda/(a\col b\col c)$ 
are representatives of
the 3-dimensional local $\overline \Lambda$-modules. According to 
Proposition 1, all the 3-dimensional $\Lambda$-modules are $\overline\Lambda$-modules,
thus the modules $M(a,b,c)$ are representatives of
the 3-dimensional local $\Lambda$-modules. 
	\medskip
{\bf (1.5)} The following theorem 
characterizes the modules of dimension at most 3 which have some relevant properties.
We write $o(q)$ for the multiplicative order of $q$.
	\medskip
{\bf  Theorem.} {\it An indecomposable module $M$ of dimension at most $3$ is 
\item{$\bullet$} torsionless if and only if  $M$ is simple or isomorphic to $\Lambda(x-y)$, to
  $\Lambda z$, to a module 
 $M(1,b,c)$ with $b\neq -q$, to $M(0,1,0)$ or to $M(0,0,1)$;
\item{$\bullet$} extensionless if and only if  $M$ is isomorphic to a module 
  $M(1,b,c)$ with $b\neq -1$;
\item{$\bullet$} reflexive if and only if  $M$ is isomorphic to a module 
  $M(1,b,c)$ with $b\neq -q^i$ for $i=1,2$;
\item{$\bullet$} Gorenstein-projective if and only if  $M$ is isomorphic to a module 
   $M(1,b,c)$ with $b\neq -q^i$
   for $i\in \Bbb Z$;
\item{$\bullet$} semi-Gorenstein-projective if and only if  $M$ is isomorphic to 
   a module  $M(1,b,c)$ with $b\neq -q^i$
   for $i\le 0$;
\item{$\bullet$} $\infty$-torsionfree if and only if  $M$ is isomorphic to 
    a module  $M(1,b,c)$
    with $b\neq -q^i$ for $i \ge 1$;
\item{$\bullet$} pivotal semi-Gorenstein-projective if and only if  $o(q) = \infty$ and 
 $M$ is isomorphic to a module  $M(1,-q,c)$;
\item{$\bullet$} pivotal $\infty$-torsionfree if and only if  $o(q) = \infty$ and 
$M$ is isomorphic to a module  $M(1,-1,c)$. 
\par}
	\medskip
 For the proof of the Theorem, see (7.9). Looking at the Theorem, the reader will be aware
that in the context considered here, 
the relevant modules of dimension at least 3 are the modules $M(1,b,c)$ with
$b,c\in k.$ Nearly all the modules mentioned in Theorem are of this kind, the only exceptions
are four isomorphism classes of torsionless modules, namely
$\Lambda(x-y)$, $\Lambda z$,  $M(0,1,0)$ and $M(0,0,1)$.

	\medskip 
{\bf (1.6)} As we have seen in (1.4), 
the set of isomorphism classes of the 3-dimensional local modules 
can be identified in a natural way 
with the projective plane $\Bbb P^2 = \Bbb P(\rad\Lambda/\rad^2\Lambda)$, with
the element $(a\col b\col c)\in \Bbb P^2$ corresponding to the module $M(a,b,c)$. 

We use homogeneous coordinates in order 
to highlight elements and subsets of $\Bbb P^2$ (or the corresponding modules):
$$
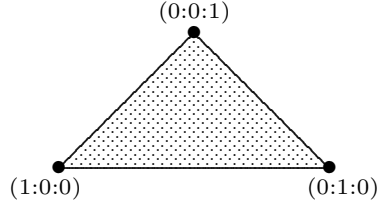
{\beginpicture
    \setcoordinatesystem units <1.5cm,1.5cm>
\multiput{$\bullet$} at 0 0  1 1  2 0 /
\plot 0 0  2 0  1 1  0 0 /
\put{$\ssize (1:0:0)$} at -.1 -.15
\put{$\ssize (0:1:0)$} at 2.1 -.15
\put{$\ssize (0:0:1)$} at 1 1.15
\setshadegrid span <.5mm>
\vshade 0 0 0  <,z,,> 1 0 1 <z,,,> 2 0 0 /
\endpicture}
$$

As Theorem (1.5) shows, of special interest is the affine subspace $H$ of $\Bbb P^2$ 
given by the points $(1\col b\col c)$ with $b,c\in k$.
As we will see in section 7, 
$H$ is a union of $\Omega\mho$-components, and
the set of 3-dimensional Gorenstein-projective 
modules is always a (proper) subset of $H$.
A module $M$ in $H$ is torsionless if and only if  it does not belong to 
the line $T = \{(1\col (-q)\col c)\mid  c\in k\}$, and is
extensionless
if and only if  it does not belong to the line $E = \{(1\col (-1)\col c)\mid c\in k\}$ (see
(6.1) and (5.1), respectively):
$$
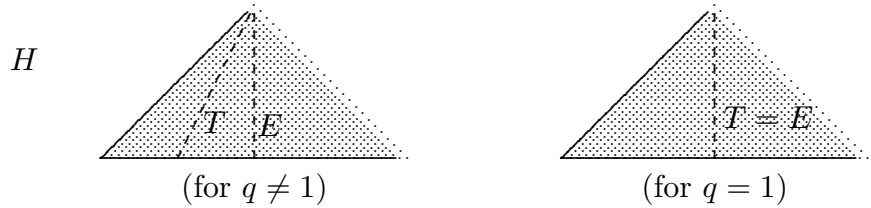
{\beginpicture
    \setcoordinatesystem units <1.7cm,1.7cm>
\put{$H$} at -1.5 0.3
\put{\beginpicture
\multiput{} at 0 0  2 1 /
\plot 1.92 0  0 0  0.95 0.95 /
\setdashes <.95mm>
\plot 1 0  1 0.95 /
\plot 0.5 0  0.975 0.95 /
\setshadegrid span <.4mm>
\vshade 0 0 0  <,z,,> 0.97 0 .97 <z,,,> 1.96 0 0 /
\put{$E$\strut} at 1.1 0.2
\put{$T$\strut} at 0.74 0.215
\setdots <1mm>
\plot  0 0  1 1  2 0  0 0 /
\put{(for $q\neq 1)$} at 1 -.2 
\endpicture} at 0 0
\put{\beginpicture
\multiput{} at 0 0  2 1 /
\plot 1.92 0  0 0  0.95 0.95 /
\setdashes <.95mm>
\plot 1 0  1 0.95 /
\setshadegrid span <.4mm>
\vshade 0 0 0  <,z,,> 0.97 0 .97 <z,,,> 1.96 0 0 /
\put{$T= E $\strut} at 1.35 0.25
\setdots <1mm>
\plot  0 0  1 1  2 0  0 0 /
\put{(for $q = 1)$} at 1 -.2 
\endpicture} at 3 0
\endpicture}
$$

{\it In case the multiplicative order $o(q)$ of $q$ is infinite,  $H$ is
the set of the $3$-dimensional modules which are
semi-Gorenstein-projective or $\infty$-torsionfree;
the line $E$ 
consists of the pivotal semi-Gorenstein-projective modules in $H$;
the line $T$ of the 
pivotal $\infty$-torsionfree modules in $H$.}

Let us emphasize: 
{\it There are $3$-dimensional pivotal semi-Gorenstein-projective modules if 
and only if  there
are $3$-dimensional pivotal 
$\infty$-torsionfree modules if and only if  the multiplicative order
of $q$ is infinite.}

	\medskip
{\bf (1.7)} The algebra $\Lambda = 
\Lambda(q)$ was exhibited in Part I in order to present in case $o(q) = \infty$
a module $M$  which is
not torsionless, such that $M$ and its $\Lambda$-dual 
$M^*$ both are semi-Gorenstein-projective:
namely the module $M = M(1,-q,0)$ with $M^* = M'(1,-q,0)$. Now we see:
	\medskip
{\it Let $o(q) = \infty$ and assume that 
$M$ is a module of dimension at most $3$. Then both $M$ and $M^*$ are
semi-Gorenstein-projective, whereas $M$ is not reflexive, if and only if $M$
is isomorphic to a module of the form $M(1,-q,c)$ with $c\in k$. In this case
$M$ is not even torsionless and all the modules $M(1,-q,c)^*$ with $c\in k$
are isomorphic.}
	\medskip
Thus, we encounter a 1-parameter family 
of pairwise non-isomorphic semi-Gorenstein-projective left modules $M$
such that their $\Lambda$-dual modules $M^*$ are isomorphic and 
semi-Gorenstein-projective, see (9.5). 
	\bigskip 

{\bf (1.8)} 
The modules $M(1,b,0)$ with $\alpha \in k$ have been studied already in Part I 
(there, they have been denoted by $M(-b)$). Theorem (1.5) shows that
these modules are quite typical for the behavior of the modules $M(1,b,c)$. 
Namely:
{\it The module $M(1,b,c)$ is Gorenstein-projective (or semi-Gorenstein-projective,
or $\infty$-torsionfree, or torsionless, or extensionless) if and only
if $M(1,b,0)$ has this property.} 
	\bigskip
{\bf (1.8) Outline of the paper.} 
Section 2 provides some preliminary results. Here, the main
target is to show that any module of length at most 3 has Loewy length
at most 2. In section 3 we collect some formulae which show that certain products
of elements in $\Lambda$ are zero. Sections 4 to 7 deal with the 
3-dimensional local left $\Lambda$-modules, section 8 with the 3-dimensional
local right $\Lambda$-modules. Section 9 discusses the $\Lambda$-duality.
The final section 10 provides an outline of the general frame for this investigation:
the study of semi-Gorenstein-projective and $\infty$-torsionfree modules over
local algebras with radical cube zero. There is an appendix which 
provides a diagrammatic description of the 3-dimensional indecomposable left 
$\Lambda$-modules. 
	\bigskip\bigskip
{\bf 2. Some left ideals and some right ideals of $\Lambda$.}
		\medskip
{\bf (2.1) Lemma.} {\it The left ideal $\Lambda (a,b,c)$ is $2$-dimensional 
if and only if  $a+b= 0$ and $ac = 0.$ We have
$\soc \Lambda(1,-1,0) = \Lambda yx$ and $\soc\Lambda(0,0,1) = \Lambda zx$.}
	\medskip
Proof. An easy calculation shows that $\soc \Lambda(1,-1,0) = \Lambda yx$ and 
$\soc\Lambda(0,0,1) = \Lambda zx$. Thus, the left ideals 
$\Lambda(0,0,1)$ and $\Lambda(1,-1,0)$ are $2$-dimensional.

Now, let $L = \Lambda(a,b,c)$ be any left ideal. If $a\neq 0$, then $yx \in L$ since
$y(a,b,c) = ayx$.

First,  assume that $a+b\neq 0.$ Then $z(a,b,c) = (a+b)zx$ shows that $zx \in L$.
We know already that for $a\neq 0,$ also $yx \in L$. If $a = 0,$ then $b\neq 0$.
Thus $x(a,b,c) = -qbyx+czx$ shows that also in this case $yx \in L.$ Thus 
$L$ cannot be 2-dimensional. 

Next, assume that $ac\neq 0.$ Since  $a\neq 0,$ we know that $yx \in L$. Since $c\neq 0,$
we use $x(a,b,c) = -qbyx+czx$ in order to see that $zx \in L.$ Again, $L$ cannot be 2-dimensional. 
$\s$
	\medskip 
{\bf (2.2)} {\it Let $L$ be a $2$-dimensional left ideal, different from
$\soc \Lambda$. Then either $L \subseteq U(1,-1,0)$ and then $\soc L = \Lambda yx$
and $L$ is isomorphic to $\Lambda(x-y)$ 
or else 
$L \subseteq U(0,0,1)$ and then $\soc L = \Lambda zx$ and $L$ is isomorphic to $\Lambda z$.}
	\medskip
Proof: There is an element $(a,b,c)+w\in L$, with $(a,b,c)\neq 0$ 
and $w\in\soc\Lambda$. Since $\rad \Lambda((a,b,c)+w) = \rad\Lambda(a,b,c)$,
also $L' = \Lambda (a,b,c)$ is $2$-dimensional and $L \subseteq L'+\soc \Lambda
= U(a,b,c)$. According to (2,1), $(a\col b\col c)$ is equal to
$(1\col(-1)\col 0)$ or to $(0\col0\col1).$ Of course, $L$ and $L'$ are isomorphic as
(left) modules. 
$\s$
	\medskip
{\bf (2.3) Lemma.} {\it There is no $3$-dimensional torsionless module with simple socle.}
	\medskip
Proof. Assume that $U$ is a $3$-dimensional torsionless module with simple socle. Then $U$
is a submodule of $\Lambda$. It is a proper submodule, thus 
of Loewy length at most 2. Therefore, 
$U$ is the sum of two 2-dimensional left ideals $L\neq L'$ with $\soc L = \soc L'.$ 
Now we use (2.2). If $L, L'$ have socle equal to $\Lambda yx$, 
then $U = L+L'= U(1,-1,0)$.
If $L, L'$ have socle equal to $\Lambda zx$, then also $U = L+L'= U(0,0,1)$. 
In both cases $\soc\Lambda \subseteq U,$ a contradiction.
$\s$
	\medskip
{\bf (2.4)} {\it Any $3$-dimensional left ideal contains $\soc \Lambda$.} 
	\medskip
{\bf (2.5)} {\it The $3$-dimensional left ideals are the subspaces 
$U(a,b,c)$. They have the following structure: $U(1,-1,0) = \Lambda(1,-1,0)\oplus
\Lambda zx$;  $U(0,0,1) = \Lambda(0,0,1)\oplus\Lambda yx$; and 
if $a+b\neq 0$ or $ac\neq 0$, then $U(a,b,c) = \Lambda(a,b,c)$ is a local
module (in particular, indecomposable).}
	\medskip
Proof. The left ideals $U(a,b,c)$ are 3-dimensional. 
Conversely, let $U$ be a 3-dimensional left ideal of $\Lambda$. 
Since $\soc \Lambda$ is contained in $U$, 
there is an element $(a,b,c)\neq 0$ with $(a,b,c)\in U$, thus $U = U(a,b,c)$. 

If $a+b= 0$ and $ac=0$, then $(a\col b\col c)$ is equal to
$(1\col (-1)\col 0)$ or to $(0\col 0\col 1)$. 
By (2.1), we have $U(1,-1,0) = \Lambda(1,-1,0)\oplus
\Lambda zx$ and $U(0,0,1) = \Lambda(0,0,1)\oplus\Lambda yx.$ 
If $a+b \neq 0$ or $ac\neq 0$, then $U(a,b,c) = \Lambda(a,b,c)$ is a local
module, thus indecomposable. 
$\s$´
	\medskip
{\bf (2.6) Proposition.} {\it Any module of dimension at most $3$ has Loewy length at most $2$.}
	\medskip
Proof. Let $M$ be a module of dimension at most 3. If $M$ is not local, then clearly
$M$ has Loewy length at most 2. If $\dim M \le 2$, then $M$ is of course 
local.  
Thus, we can assume that $M$ is 3-dimensional and local and therefore a factor module
of $\Lambda$, say $M = \Lambda/U$. According to (2.4), $\soc \Lambda \subseteq U$,
thus $M$ is annihilated by
$\soc\Lambda$, and therefore $M$ has Loewy length at most 2.
$\s$
	\medskip
{\bf (2.7) Lemma.} {\it Any indecomposable torsionless module $M$ of dimension at most $3$
is local and isomorphic to a left ideal of $\Lambda$. If $\dim M = 3$, then $M$ 
is of the form $U(a,b,c)$.}
	\medskip
Proof. Let $M$ be indecomposable and torsionless. If $\dim M \le 2$, then $M$ is of course 
local and isomorphic to a left ideal. Thus we can assume that $\dim M = 3.$

Since $M$ is torsionless, there is a set of non-zero 
maps $u_i\:M \to {}_\Lambda\Lambda$ (say with index set $I$)
such that $\bigcap_{i\in I} K_i = 0$, where $K_i$ is the kernel of $u_i.$ 

If $K_i = 0$ for some $i$, then already $u_i$ is an embedding (thus $M$ is
isomorphic to a left ideal). 
In particular, if the socle of $M$ is simple, then we must have $K_i = 0$ for some $i$.
Thus, we can assume that the socle of $M$ is not simple. Therefore $M$ has to be
a local module and we have a surjective map $\pi\:{}_\Lambda\Lambda \to M$.

It remains to look at the case where $\dim K_i = 1$ or $2$ for all $i$. Since
the only 2-dimensional submodule of $M$ is its radical, we have 
$\bigcap_{i\in I'} K_i = 0,$ where $I'$ is the set of indices $i$ with
$\dim K_i = 1$. But then $K_i\cap K_j = 0$ for some $i\neq j$ in $I'$. 
This shows that we can assume that  $I = \{1,2\}$ and that $K_1, K_2$
are different 1-dimensional submodules of $M$. 

Now $u_i$ provides an isomorphism from $M/K_i$ onto 
a (2-dimensional) left ideal of $\Lambda$. Since $M/K_i$ is 
indecomposable, (2.2) shows that $M/K_i$ is isomorphic to 
$\Lambda(1,-1,0)$ or to $\Lambda(0,0,1)$. Let $K'_i = \Ker(u_i\pi)$
for $i=1,2$.

If $M/K_i \simeq \Lambda(1,-1,0)$, then $K'_i$ is equal to
$\Lambda(x+qy)+\Lambda z$, since $\Lambda(0,0,1)$ is annihilated by
$x+qy$ and by $z$. Similarly, 
if $M/K_i \simeq \Lambda(0,0,1)$, then $K'_i$  is equal  to  
$\Lambda(x+qy)+\Lambda z$. Thus one of $M/K_i$ has to be isomorphic
to $\Lambda(1,-1,0)$, the other one to $\Lambda(0,0,1)$ and
$\Ker(\pi) = K_1'\cap K'_2 = U(0,0,1).$ It follows that
$M \simeq {}_\Lambda\Lambda/\Ker(\pi) = {}_\Lambda\Lambda/U(0,0,1).$
But ${}_\Lambda\Lambda/U(0,0,1)$ is isomorphic to 
the left ideal $\Lambda(1,-1,1) = U(1,-1,1)$.

We have shown that $M$ is isomorphic to a left ideal, thus of the form $U(a,b,c)$,
see (2.5). Since we assume that
$M$ is indecomposable, (2.5) asserts that $M$ is local.
$\s$
	\bigskip

We need to know also the right ideals $(a,b,c)\Lambda$. Note that $U(a,b,c)$
is always a twosided ideal and it will be pertinent to denote $U(a,b,c)$ by $U'(a,b,c)$, 
if we consider it as a right ideal (thus as a right module). 
	\medskip
{\bf (2.8) The right ideals $(a,b,c)\Lambda$}. 
{\it If $a \neq 0$ or $bc \neq 0$, then $(a,b,c)\Lambda = U(a,b,c)$ is $3$-dimensional.
The right ideals  $(0,1,0)\Lambda$ and $(0,0,1)\Lambda$ are  $2$-dimensional with
$\soc\, (0,1,0)\Lambda = yx\Lambda$
and $\soc\, (0,0,1)\Lambda = zx\Lambda.$}
	\medskip
Proof: Let $V = (a,b,c)\Lambda$. First, let
$a\neq 0$. Then $zx$ belongs to $V$, since $(a,b,c)z = azx.$ Also $yx \in V$, since
$(a,b,c)y = -qayx + czx$. Second, assume that $a = 0$ and $bc \neq 0.$ Then 
$(0,b,c)y = czx$ shows that $zx \in V$, and $(0,b,c)x = byx + czx$ shows that also $yx \in V.$
$\s$
	\medskip
{\bf (2.9)} {\it If a $3$-dimensional indecomposable right module $N$
is torsionless, then it is isomorphic to a right ideal, thus to $U'(a,b,c)$ for some
$(a,b,c)\neq 0.$} 
	\medskip
Proof. Let $N$ be a 3-dimensional indecomposable torsionless right module.
As in (2.7) one shows that $N$ is isomorphic to a right ideal, 
using (2.8) instead of (2.2). It remains to see that all 3-dimensional right ideals
are of the form $U'(a,b,c)$. Here, one has to copy the proof of (2.5).
	\bigskip\bigskip
{\bf 3. The transformations $\omega$ and $\omega'$.}
	\medskip
If $(a\col b\col c)$ is different from $(1\col (-1)\col 0)$ and $(0\col 0\col 1)$, then
(2.5) shows that  $U(a,b,c)$ is a $3$-dimensional local module,
thus of the form $M(a'\col b'\col c')$. In order to describe in which way 
$(a'\col b'\col c')$ depends on $(a\col b\col c)$, we will need the 
transformations $\omega$ and $\omega'$.
We start with some equalities in $\Lambda$.
	\medskip
{\bf (3.1)} {\bf Formulae.} {\it Let $a,b,c\in k.$ Then
$$
\alignat 2
 \bigl(a\,x+qb\,y-\tfrac a{a+b}c\,z\bigr)\bigl(a\,x + b\,y + c\,z\bigr) &= 0  &
  \quad\text{if}\quad & a+b\neq 0  
\tag1 \cr
 z(ax-ay+cz) &=0 && \tag2 \cr
 \bigl(a\,x + b\,y + c\,z\bigr)\bigl(a\,x+q^{-1}b\,y-\tfrac {a+q^{-1}b}{a}c\,z\bigr) &= 0
 &
  \quad\text{if}\quad & a\neq 0  
\tag3 \cr
 (by+cz)z &=0 && \tag4
\endalignat 
$$}

Proof of the equality (1):
$$
\align
  \bigl(a\,x+&qb\,y-\tfrac a{a+b}c\,z\bigr)\bigl(a\,x + b\,y + c\,z\bigr) \cr
 &=\ ab\,xy+ac\,xz+qab\,yx-\tfrac a{a+b}ac\,zx - \tfrac{a}{a+b}bc\,zy \cr
 &=\ ab\,(xy+q\,yx) +\bigl(1- \tfrac a{a+b}-\tfrac{b}{a+b}\bigr)ac\,zx \ =\  0.
\endalign
$$
The proof of the remaining equalities is similar.  $\s$
	\bigskip
{\bf (3.2)} In case $a+b\neq 0$, let 
$ \omega(a,b,c) = (a,qb,-\tfrac a{a+b}c).$
In case $a'\neq 0,$ let $\omega'(a',b',c') = (a',q^{-1}b', -\tfrac {a'+q^{-1}b'}{a'}c').$
	\bigskip
{\bf Proposition.} {\it The transformation  
$\omega$ provides a bijection from the set $\{(a,b,c)\in k^3 \mid a(a+b)\neq 0\}$ 
onto the set $\{(a',b',c')\in k^3 \mid a'(a'+q^{-1}b')\neq 0\},$ with inverse $\omega'$.} 
	\medskip
Proof. Let $a(a+b)\neq 0$. Then $(a',b',c') = \omega(a,b,c)$ is defined and
$a' = a \neq 0$, and $a'+q^{-1}b' = a+q^{-1}qb = a+b \neq 0.$ Thus
$\omega$ maps  $\{(a,b,c)\in k^3 \mid a(a+b)\neq 0\}$ into 
$\{(a,b,c)\in k^3 \mid a(a+q^{-1}b)\neq 0.$ Similarly, $\omega'$ maps
$\{(a',b',c')\in k^3 \mid a'(a'+q^{-1}b')\neq 0$ into 
$\{(a,b,c)\in k^3 \mid a(a+b)\neq 0.$ It is easy to check that 
$\omega'\omega(a,b,c) = (a,b,c)$ for $a(a+b)\neq 0$ 
and that $\omega\omega'(a',b',c') = (a',b',c')$ for $a'(a'+q^{-1}b')\neq 0.$
$\s$
	\bigskip
{\bf 4. The isomorphism class of $U(a,b,c) \simeq \Omega M(a,b,c)$.}
	\medskip
{\bf (4.1) Proposition.} {\it Let $(a,b,c)\neq 0.$ Then 
$$
\Omega M(a,b,c) \simeq \left\{
\matrix M(\omega(a,b,c))                &&\text{if}  
                                         & a\neq 0,\ a+b\neq 0,\qquad\quad &\qquad & (1)\cr
        M(0,0,1)                        &&\text{if}  
                                         & a\neq 0,\  a+b =0,\ c\neq 0,  &&(2)\cr
        \Lambda (x-y)\oplus \Lambda zx  &&\text{if}  
                                         & a\neq 0,\ a+b = 0,\ c = 0, &&(3)\cr  
        M(0,1,0)                        &&\text{if}  
                                         & \ a = 0,\ b \neq 0,\qquad\qquad\quad   &&(4) \cr
        \Lambda z \oplus    \Lambda yx  &&\text{if}  
                                         & \ a = 0,\ b = 0.\qquad\qquad\quad      &&(5)\cr
\endmatrix
\right.
$$}

Proof: If $a = 0$ and $b = 0,$ then
$U(a,b,c) = U(0,0,1)$. If $a+b= 0$ and $c = 0$, then $U(a,b,c) = U(1,-1,0)$.
According to (2.3), $U(0,0,1) =  \Lambda z \oplus    \Lambda yx$ and 
$U(1,-1,0) = \Lambda (x-y)\oplus \Lambda zx$, This shows (5) and (3). In this way, we
have considered all triples $(a,b,c)$ with $a+b = 0$ and $ac = 0$. 

Thus, let $a+b\neq 0$ or $ac \neq 0$. By (2.5), $U(a,b,c) = \Lambda(a,b,c)$ is
local and we look at the surjective 
map $\phi\:{}_\Lambda\Lambda \to U(a,b,c)$ which sends $1$ to $(a,b,c)$.

Let $a+b\neq 0.$ According to formula (1) of (3.1), $\Lambda(a,b,c)$ is annihilated
by $\omega(a,b,c)$, thus $M(\omega(a,b,c))) = {}_\Lambda\Lambda/\Lambda(\omega(a,b,c))$
maps onto $\Lambda(a,b,c)$. Since the modules $M(\omega(a,b,c))$ and
$\Lambda(a,b,c)$ both have dimension 3, we see that 
$U(a,b,c) = \Lambda(a,b,c)$ is isomorphic to $M(\omega(a,b,c))$. This yields (1)
and (4) (namely, if $a = 0,$ and $b\neq 0$, we have $\omega(0,b,c) = (0,qb,0)$). 

Finally, we show (2). For $c\neq 0$, the module $U(1,-1,c)$ is isomorphic
to $M(0,0,1)$. Now we use in the same way formula (2) of (3.1). 
$\s$
	\medskip
The following picture outlines the position of the partition of $\Bbb P^2$
which is used in the Proposition.
$$
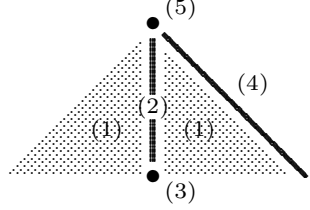
{\beginpicture
    \setcoordinatesystem units <1.7cm,1.7cm>
\multiput{} at 0 0  2 1 /
\multiput{$\bullet$} at 1 1  1 0 /
\multiput{$\ssize(1)$} at  0.7 0.3  1.3 0.3 /
\plot 1 0.1  1 0.37 /
\plot 0.985 0.1  0.985 0.37 /
\plot 1.01 0.1  1.01 0.37 /
\plot 1 0.54  1 0.9 /
\plot 0.985 0.54  0.985 0.9 /
\plot 1.01 0.54  1.01 0.9 /
  
\plot 1.07 0.93  2 0 /
\plot 1.077 0.937   2.007 0.007 /
\plot 1.063 0.923  1.993 -0.007 /
\setshadegrid span <.4mm>
\vshade 0 0 0  0.95 0 0.95 /
\vshade 1.05 0 0.85  1.93 0 0 /
\put{$\ssize (5)$} at 1.18 1.1
\put{$\ssize (4)$} at 1.65 0.6
\put{$\ssize (2)$} at 1 0.45
\put{$\ssize (3)$} at 1.18 -.1
\endpicture}
$$

	\medskip
{\bf (4.2) Corollary.} {\it The syzygy functor 
$\Omega$ provides a bijection from the set 
of isomorphism classes of modules $M(a,b,c)$ with $a(a+b)\neq 0$ onto 
the set of isomorphism classes of modules $M(a',b',c')$ with $a'(a'+q^{-1}b')\neq 0$
and we have $\Omega M(a,b,c) = M(\omega(a,b,c))$ for $a(a+b)\neq 0$.}
	\medskip
Proof. This follows directly from Propositions (3.2) and (4.1). 
$\s$
	\bigskip

{\bf 5. The extensionless modules $M(a,b,c)$.}
	\medskip
{\bf (5.1) Proposition.} {\it The module $M(a,b,c)$ is extensionless if and only if  $a(a+b)\neq 0.$}
	\medskip
For the proof, we need some preparations.
	\medskip
{\bf (5.2) Lemma.} {\it The following conditions are equivalent:
\item{\rm (i)} The module $M(a,b,c)$ is extensionless.
\item{\rm(ii)} The inclusion map $\iota\:U(a,b,c) \to {}_\Lambda\Lambda$ is a left
 $\add(\Lambda)$-approximation.
\item{\rm(iii)} $U(a,b,c) = \Lambda(a,b,c)$ and 
 the inclusion map $\iota\:\Lambda(a,b,c) \to {}_\Lambda\Lambda$ is a left
 $\add(\Lambda)$-approximation.
\item{\rm(iv)} The subspace $U(a,b,c)$ is 
  indecomposable both as a left module and as a right module,
  and the image of every homomorphism ${}_\Lambda U(a,b,c) \to {}_\Lambda\Lambda$ is 
  contained in $U(a,b,c)$.\par}
	\medskip
Proof. The equivalence of (i) and (ii) follows from Part I, Lemma 2.1.

(ii) $\implies$ (iii): We assume (ii).
If $U(a,b,c) = U_1\oplus U_2$ with $U_1,U_2$ both non-zero, then a minimal left 
$\add(\Lambda)$-approximation $U(a,b,c)\to \Lambda^t$ is the direct sum of
minimal left $\add(\Lambda)$-approximations $U_1 \to \Lambda^{t_1}$ and  
$U_2 \to \Lambda^{t_2}$, thus $t = t_1+t_2 \ge 2.$ This shows that $U(a,b,c)$ is
indecomposable. According to (2.5), this means that $U(a,b,c) = \Lambda(a,b,c)$.

(iii) $\implies$ (iv). Since $\Lambda(a,b,c)$ is a local module, it is indecomposable.
Thus $U(a,b,c) = \Lambda(a,b,c)$ implies that $U(a,b,c)$ considered as a left module
is indecomposable. Given any homomorphism $\phi\:U(a,b,c) \to {}_\Lambda\Lambda$, (iii)
provides $\lambda\in \Lambda$ with $\phi(a,b,c) = (a,b,c)\lambda \in (a,b,c)\Lambda 
\subseteq U(a,b,c).$ Now assume that $(a,b,c)\Lambda$ is a proper subset of
$U(a,b,c).$ Let $w\in \soc\Lambda$. Since $\Lambda w$ is simple, there is 
a homomorphism $\phi\:\Lambda(a,b,c) \to \Lambda$ with $\phi(a,b,c) = w$
and (iii) asserts that $w = \phi(a,b,c) = (a,b,c)\lambda$ for some $\lambda\in
\Lambda$. This shows that $\soc \Lambda \subseteq (a,b,c)\Lambda$ and therefore
$U(a,b,c) = (a,b,c)\Lambda$. In particular, $U(a,b,c)$ is indecomposable also as a right
$\Lambda$-module.

(iv) $\implies$ (ii). Let $\phi\:U(a,b,c) \to {}_\Lambda\Lambda$ be a homomorphism.
Since $U(a,b,c)$ is indecomposable as a left module, we have $U(a,b,c) = \Lambda(a,b,c)$.
Since $U(a,b,c)$ is indecomposable as a right module, we have $U(a,b,c) = (a,b,c)\Lambda$.
According to (iv), $\phi(a,b,c)\in U(a,b,c) = (a,b,c)\Lambda,$ thus $\phi(a,b,c) =
(a,b,c)\lambda = r_\lambda\iota(a,b,c)$ for some $\lambda\in \Lambda$, where
$r_\lambda\:{}_\Lambda\Lambda \to {}_\Lambda\Lambda $ is the right multiplication by $\lambda$.
Since the left module $U(a,b,c) = \Lambda(a,b,c)$ is generated by $(a,b,c)$, the
equality $\phi(a,b,c) = r_\lambda\iota(a,b,c)$ implies that $\phi = r_\lambda\iota.$ 
$\s$
	\medskip
{\bf (5.3) Lemma.} {\it Let $R$ be a ring and $X$ a left $R$-module.
If $\phi\:{}_RR \to X$ is an $R$-module 
homomorphism and $w\in R$
annihilates $X$, then $R w \subseteq \Ker\phi.$}  
	\medskip
{\bf Corollary.} {\it Let $L$ be a left ideal of $R$ and $X$ 
an $R$-module annihilated by $w_1,\dots,w_t\in R$.
The image of any map $R/L \to X$ is a factor module
of $R/(L+R w_1+\cdots R w_t)$.}
	\medskip
Proof. Let $\phi\:R/L \to X$ be a homomorphism.
Let $\pi\:R \to R/L$ be the canonical projection. 
By construction, $L$ is contained in
$\Ker(\phi\pi)$. By the lemma, also the left ideals $R w_i$ are contained in
$\Ker(\phi\pi).$ Thus $L+R w_1+\cdots +R w_t \subseteq \Ker(\phi\pi).$
$\s$
	\bigskip
{\bf (5.4)} Proof of Proposition (5.1).
According to (5.2), $M(a,b,c)$ is extensionless if and only if  condition (iv)
is satisfied. We look at all the elements $(a\col b \col c) \in \Bbb P^2,$
using the partition of $\Bbb P^2$ into the subsets (1) to (5) as in (4.1).

The cases (3) and (5): Both $U(1,-1,0)$ and $U(0,0,1)$ are decomposable 
as left modules, see (2.5).
Case (4): According to (4.1), $U(0,1,c)  
\simeq M(0,1,0)$. Obviously, $M(0,1,0)$ has $\Lambda z$ as a factor module, thus
there is a homomorphism $U(0,1,c) \to
{}_\Lambda\Lambda$ with image $\Lambda z$ and $\Lambda z \not\subseteq U(0,1,c)$.
The case (2) is similar: (4.1) shows that $U(1,-1,c) \simeq M(0,0,1)$, 
and $M(0,0,1)$ maps onto
$\Lambda z$; thus there is a homomorphism $U(1,-1,c) \to
{}_\Lambda\Lambda$ with image $\Lambda z$ and $\Lambda z \not\subseteq U(1,-1,c)$.
This shows that none of the modules $M(a,b,c)$ with $a(a+b) = 0$ is extensionless. 
	\smallskip
It remains to consider the case (1). Thus, assume that $a(a+b)\neq 0.$ 
Let $(1,b',c') = \omega(1,b,c)$, thus $b' = qb.$
We want to show that the conditions (iv)
of (5.2) are satisfied. According to (2.5) and (2.8), 
$U(a,b,c)$ is indecomposable both as a left module and as a right module,
It remains to show that 
the image of every homomorphism ${}_\Lambda U(a,b,c) \to {}_\Lambda\Lambda$ is 
contained in $U(a,b,c)$.
	\smallskip
(a) {\it The only left ideal isomorphic to $U(1,b,c)$ is $U(1,b,c)$
itself.}
Proof. The 3-dimensional left ideals are of the form $U(a'',b'',c'')$, for some
$(a'',b'',c'')\neq 0$, see (2.5).
Assume that $U(1,b,c) \simeq U(a'',b'',c'')$. We have $U(a'',b'',c'') \simeq
\Omega M(a'',b'',c'')$ and by (4.1) we must be in case (1), namely
$a''\neq 0$ and $a''+b''\neq 0$. In particular, we may assume that
$a'' = 1$ and (4.1)(1) asserts that  
$\Omega M(1,b'',c'') = 
M(\omega(1,b'',c'')).$ The isomorphy $M(\omega(1,b,c)) \simeq M(\omega(1,b'',c''))$
implies that the triples $\omega(1,b,c)$ and $\omega(1,b'',c'')$ yield the same
element in $\Bbb P^2,$ and since the first coordinate of both triples is equal to 1,
we have $\omega(1,b,c) = \omega(1,b'',c'').$ Since $1+b \neq 0$ and $1+b''\neq 0,$
we use (3.2) in oder to conclude that $(1,b,c) = (1,b'',c'').$
	\smallskip
(b) {\it The left ideal $\Lambda z$ is not a factor module of $U(1,b,c)$.}
The proof uses Corollary (5.3) for the left ideal $L = U(1,b',c')$ and the module
$X = \Lambda z$ which is annihilated by $y$ and $z$. 
Namely, on the one hand, we have 
$U(1,b,c) \simeq \Omega M(1,b,c) \simeq M(\omega(1, b, c)) = M(1,b',c') = 
\Lambda/U(1, b', c') = \Lambda/L.$
On the other hand, $\rad \Lambda = \Lambda(x+b'y+c'z)+\Lambda y +\Lambda z \subseteq 
U(1,b',c')+\Lambda y +\Lambda z \subseteq \rad\Lambda$ shows that 
$L+\Lambda y + \Lambda z = \rad \Lambda$. Therefore, (5.3) asserts that the image 
of any homomorphism $U(1,b,c) \to \Lambda z$
is a factor module of $\Lambda/\rad\Lambda$, thus simple or zero.
	\smallskip
(c) {\it The left ideal $\Lambda(x-y)$ is not a factor module of $U(1,b,c)$.}
Again, we use Corollary (5.3) for $L = U(1,b',c')$ and now for 
$X = \Lambda (x-y)$. Note that $\Lambda (x-y)$ is annihilated by $x-qy$ and $z$.
We recall from (b) that $U(1,b,c) \simeq \Lambda/L.$
And we have $\rad \Lambda = \Lambda(x+b'y+c'z)+\Lambda(x-qy) +\Lambda z$,
since $b' = qb \neq -q$. Therefore, we also have
$U(1,b',c')+\Lambda(x-qy) +\Lambda z = \rad\Lambda$, and 
(5.3) asserts that the image of any homomorphism $U(1,b,c) \to \Lambda z$
is simple or zero.
	\smallskip 
{\it Any homomorphism $\phi\:U(1,b,c) \to {}_\Lambda\Lambda$ maps into $U(1,b,c)$.}
Proof. According to (b) and (c), the image $I$ of $\phi$ 
is not of dimension $2$. If the image $I$ is 
of dimension $3$, then (a) shows that $I$ is equal to $U(1,b,c)$. Of course, if $I$ 
is of dimension at most $1$, then $I \subseteq \soc \Lambda \subseteq U(1,b,c).$
$\s$
	\medskip
{\bf (5.5) Corollary.} {\it If $M(a,b,c)$ is extensionless, then $\Omega M(a,b,c) 
\simeq M(\omega(a,b,c))$.}
	\medskip
Proof. This follows directly from (5.1) and the case (1) of (4.1).
$\s$.
	\bigskip

{\bf 6. The torsionless modules $M(a,b,c)$.}
	\medskip
{\bf (6.1) Proposition.} {\it The module $M(a,b,c)$ is torsionless if and only if  either
   $a(a+q^{-1}b)\neq 0$ or else $a = 0$ and $bc = 0$} (so that
$(a\col b \col c)$ is equal to $(0\col 1\col 0)$
or to $(0\col 0 \col 1)$).
	\medskip
In order to prove (6.1), we consider the possible cases separately.
First, we consider the modules $M(a,b,c)$ with $a\neq 0.$
In section 5 we have seen that $M(1,b,c)$ is extensionless if and only if  $b\neq -1$, and 
then $\Omega M(1,b,c) \simeq M(\omega(1,b,c))$. There is the following corresponding assertion 
concerning the torsionless modules (see also (7.1)).
	\medskip
{\bf (6.2)} 
{\it  The module $M(1,b,c)$ is torsionless if and only if  $b\neq -q$, and in this case
$\mho M(1,b,c) \simeq M(\omega'(1,b,c))$.}
	\medskip
Proof. Let $b\neq -q$. Then $\omega'(1,b,c) = (1,q^{-1}b,c')$ for some $c'$.
According to (5.1) and (5.5), $M(1,q^{-1}b,c')$ is extensionless and
$\Omega M(1,q^{-1}b,c') \simeq M(1,b,c)$,  
since $\omega(1,q^{-1}b,c') = \omega\omega'(1,b,c) = (1,b,c).$
This shows that $M(1,b,c)$ is torsionless and that $\mho M(1,b,c) \simeq M(\omega'(1,b,c))$.
 
Conversely, we consider $M(1,-q,c)$ and assume, for the contrary, that $M(1,-q,c)$ 
is torsionless. 
According to (2.7), this means that $M(1,-q,c)$ is isomorphic to 
a left ideal $U(a',b',c') = \Omega M(a',b',c').$ According to
(4.1), we must be in the case $a'+b' \neq 0$ and $a'\neq 0$.
We can assume that $a' = 1$, thus $1+b'\neq 0$. We have
$\Omega M(1,b',c') \simeq M(\omega(1,b',c')) = M(1,qb',c'')$ 
for some $c''$. Since $M(1,-q,c) \simeq 
\Omega M(1,b',c') \simeq M(1,qb',c''),$ we see 
that $(1,-q,c) = (1,qb',c'')$, thus $b' = -1.$ But this is
a contradiction to $1+b'\neq 0.$
$\s$
	\medskip 
{\bf (6.3)} {\it For $M = M(0,1,0)$ and $M(0,0,1)$, there is no monomorphism 
$M \to {}_\Lambda\Lambda$ which is an $\add(\Lambda)$-approximation.}
	\medskip
Proof. Let $M$ be equal to $M(0,1,0)$ or to $M(0,0,1)$.
Assume that there is a monomorphism $u\:M \to {}_\Lambda\Lambda$ which is an 
$\add(\Lambda)$-approximation. The image $u(M)$ is a 3-dimensional left ideal,
thus of the form $U(a,b,c)$ for some $(a,b,c)\neq 0$, see (2.7).
The implication (ii) $\implies$ (iv) in (5.2) asserts that any
homomorphism $U(a,b,c) \to {}_\Lambda\Lambda$ maps into $U(a,b,c)$.

Obviously, both modules $M(0,1,0)$ and $M(0,0,1)$ have a factor module 
isomorphic to $\Lambda z,$ thus there is a surjective homomorphism $U(a,b,c) \to 
\Lambda z,$ and therefore $\Lambda z \subseteq U(a,b,c)$. But $\Lambda z$
is an indecomposable module of length 2, and $U(a,b,c) \simeq M$ is a local module
of length 3
with socle of length 2. A local module of length 3 
with socle of length 2 has no indecomposable
submodule of length 2, thus we obtain a contradiction.
$\s$
	\bigskip
{\bf (6.4) Proposition.} {\it The modules $M(0,b,c)$ with $bc\neq 0$
are not torsionless.}
	\medskip
Proof. Let $M = M(0,b,c)$ with $bc\neq 0$ and assume that $M$ is torsionless.
According to (2.7), this means that $M \simeq U(a',b',c') \simeq \Omega M(a',b',c')$
for some triple $(a',b',c'),$  and (2.5) asserts that $a'+b'\neq 0$ or $a'c'\neq 0$.
Now we use (4.1) and have to distinguish the three cases (1), (2) and (4).
Case (1) means that $a'+b'\neq 0$ and $a'\neq 0$, then $\Omega M(a',b',c') \simeq 
M(\omega(a',b',c'))$ and the first component of $\omega(a',b',c')$ is $a'$, thus
non-zero. But then $M(\omega(a',b',c'))$ cannot be isomorphic to $M(0,b,c)$.
Case (4) means that  $a' = 0$ and $b' \neq 0$. Then $\Omega M(a',b',c') \simeq M(0,1,0)$, 
thus not isomorphic to $M(0,b,c)$ with $bc\neq 0.$ Finally, there is the
case (2) with $a'+b' = 0$ and $a'c' \neq 0$. Then $\Omega M(a',b',c') \simeq M(0,0,1)$, 
again not isomorphic to $M(0,b,c)$ with $bc\neq 0.$ 
In all cases, we get a contradiction. 
$\s$
	\medskip
{\bf (6.5) Proposition.} {\it  If $M$ is equal to
$M(0,1,0)$ or $M(0,0,1)$, then $M$ is torsionless and the module $\mho M$ has
Loewy length $3$. Since $\mho M$ is indecomposable and non-projective, it is not torsionless.}
	\medskip 
Proof. 
The modules $M$ of the form $M(0,1,0)$ and $M(0,0,1)$ are torsionless, since 
(4.1), (4) and (2) 
assert that $M(0,1,0) \simeq \Omega M(0,1,0)$ and that $M(0,0,1) \simeq \Omega M(1,-1,1)$.
According to (5.2), in both cases there is no inclusion map $M \to \Lambda$ which is
an $\add(\Lambda)$-approximation. Thus, a minimal left 
$\add(\Lambda)$-approximation of $M$ is an injective map $M \to \Lambda^t$ 
with $t \ge 2.$ This shows that $\mho M$ has dimension $6t-3$ and its top has dimension $t$. 
According to Part I (3.2), $\mho M$ is indecomposable and not projective. 
The Loewy length of $\mho M$ has to be 3. [Namely, 
an indecomposable module with Loewy length at most $2$ and top of dimension $t\ge 2$ 
has dimension at most $4t-1$, since it is a proper factor module of 
$\overline\Lambda^t$. But $6t-3\le 4t-1$ implies $t\le 1,$ a contradiction.]
An indecomposable non-projective module of Loewy
length 3 cannot be torsionless.  
$\s$
	\bigskip
{\bf (6.6)} We finish this section by reformulating the results concerning the modules of
the form $M(0,b,c)$ in terms of $\Omega\mho$-components. Here,
we will exhibit the structure of all the $\Omega\mho$-components containing 
modules of the form
$M(0,b,c)$.  We have to distinguish between 
the modules $M(0,1,0)$ and $M(0,0,1)$ and the modules $M(0,b,c)$ with $bc\neq 0$, thus
lying on the dashed line $A' = \{(0\col b\col c)\mid bc\neq 0\}$:
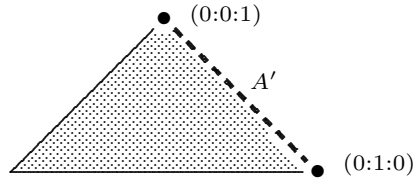
$$
{\beginpicture
    \setcoordinatesystem units <1.7cm,1.7cm>
\multiput{} at 0 0  2 1 /
\multiput{$\bullet$} at 1 1  2 0 /
\plot 0 0  0.93 0.93 /
\plot 0 0  1.9 0 /
\setdashes <1mm>
\plot 1.07 0.93  1.93 0.07 /
\plot 1.077 0.937  1.937 0.077 /
\plot 1.063 0.923  1.923 0.063 /
\setshadegrid span <.4mm>
\vshade 0 0 0  <,z,,> 0.95 0 .91 <z,z,,> 1.06 0  0.9  <z,,,> 1.92 0 0 /
\put{$\ssize (0: 0:1)$} at 1.4 1.02
\put{$\ssize A'$} at 1.65 0.6
\put{$\ssize (0:1:0)$} at 2.4 0.05
\endpicture}
$$

{\it The modules in $A'$ are singletons (that is,
components of type $\Bbb A_1$) in the $\Omega\mho$-quiver. And, there are
the following two $\Omega\mho$-components of the form $\Bbb A_2$:}
$$
{\beginpicture
    \setcoordinatesystem units <2cm,1cm>
\put{\beginpicture
\put{$\sssize \blacksquare$} at 1 0 
\put{$\ssize \blacklozenge$} at 0 0 
\arra{0.15 0}{0.1 0}
\setdashes <1mm>
\arr{0.9 0}{0.1 0}
\put{$M(0,0,1)$} at -.2 -.3
\put{$\mho M(0,0,1)$} at 1.2 -.3
\endpicture} at 0 0
\put{\beginpicture
\put{$\sssize \blacksquare$} at 1 0 
\put{$\ssize \blacklozenge$} at 0 0 
\arra{0.15 0}{0.1 0}
\setdashes <1mm>
\arr{0.9 0}{0.1 0}
\put{$M(0,1,0)$} at -.2 -.3
\put{$\mho M(0,1,0)$} at 1.2 -.3
\endpicture} at 3 0
\endpicture}
$$
(If $M$ is an indecomposable module,
then we represent $[M]$ in the $\Omega\mho$-quiver usually just by a circle $\circ$.
We use a bullet $\bullet$ in case  we know
that $M$ is torsionless and extensionless, a black square $\sssize \blacksquare$ 
in case we know that $M$
is extensionless, but not torsionless; and a black lozenge $\ssize\blacklozenge$ in case we know
that $M$ is torsionless, but not extensionless.)
	\bigskip\bigskip
{\bf 7. The modules $M(1,b,c)$ and proof of Theorem (1.5).}
	\medskip
We consider now the affine subspace $H$ of $\Bbb P^2$ 
given by the points $(1\col b\col c)$ with $b,c\in k$ and the corresponding
modules $M(1,b,c)$. We recall that $o(q)$ denotes the multiplicative order of $q$.
	\medskip
{\bf (7.1)} We have seen in (4.2) that $\Omega$ provides a bijection from the set of  modules $M(1,b,c)$ with $b\neq -1$ onto the set of modules $M(1,b',c')$ with $b'\neq -q$.
The sections 5 and 6 strengthen this bijection as follows:
	\medskip
{\it If $b\neq -1$, then the exact sequence
$$
 0 \to M(1,b',c') \to {}_\Lambda\Lambda \to M(1,b,c) \to 0
$$
with $(1,b',c') = \omega(1,b,c)$ is an $\Omega\mho$-sequences (here, $(1,b',c')$
is an arbitrary triple with $b'\neq -q,$ and $(1,b,c) = \omega'(1,b',c')$).
We obtain in this way
all the $\Omega\mho$-sequences involving modules of the form
$M(1,b,c)$.} 
	\bigskip
{\bf (7.2) Reformulation.} 
The neighborhood of $M(1,b,c)$ in the $\Omega\mho$-quiver
looks as follows:
$$
{\beginpicture
    \setcoordinatesystem units <2cm,.65cm>
\put{\beginpicture
    \setcoordinatesystem units <2cm,1.1cm>
\multiput{$\bullet$} at  1 0  /
\multiput{$\circ$} at  0 0  2 0 /
\arra{0.15 0}{0.1 0}
\arra{1.15 0}{1.1 0}
\setdashes <1mm>
\arr{0.9 0}{0.1 0}
\arr{1.9 0}{1.1 0}
\put{$M(1,b,c)$} at 1 -.3
\put{$M(\omega(1,b,c))$} at -.2 -.3
\put{$M(\omega'(1,b,c))$} at 2.2 -.3
\multiput{$\dots$} at -.35 0  2.35 0 /
\multiput{} at -1 0  3 0 /
\endpicture} at 0 0
\put{\beginpicture
    \setcoordinatesystem units <2cm,1.1cm>
\put{$\sssize \blacksquare$} at 1 0 
\put{$\circ$} at 0 0
\arra{0.15 0}{0.1 0}
\setdashes <1mm>
\arr{0.9 0}{0.1 0}
\put{$M(1,-q,c)$} at 1 -.3
\put{$M(\omega(1,-q,c))$} at -.2 -.3
\multiput{$\dots$} at -.35 0  /
\multiput{} at -1 0  3 0 /
\endpicture} at -.4 -1.5
\put{\beginpicture
    \setcoordinatesystem units <2cm,1.1cm>
\put{$\ssize \blacklozenge$} at 1 0 
\put{$\circ$} at  2 0
\arra{1.15 0}{1.1 0}
\setdashes <1mm>
\arr{1.9 0}{1.1 0}
\put{$M(1,-1,c)$} at 1 -.3
\put{$M(\omega'(1,-1,c))$} at 2.2 -.3
\multiput{$\dots$} at  2.35 0 /
\multiput{} at -1 0  3 0 /
\endpicture} at 0.4 -3

\put{$b\notin\{-1,-q\}$} at 3.1 0
\put{$b = -q \neq -1$} at 3.1 -1.5
\put{$b = -1\neq -q$} at 3.1 -3

\endpicture}
$$
and $M(1,b,c)$ is a singleton in the $\Omega\mho$-quiver if $q = 1$ and  $b=-1$.
	\medskip
{\bf (7.3)}
{\it The module
$M(1,b,c)$ is semi-Gorenstein-projective if and only if  $b\neq -q^t$
for all $t\le 0$.
The module
$M(1,b,c)$ is $\infty$-torsionfree if and only if  $b\neq -q^t$
for all $t\ge 1$.}
	\medskip
Proof: $M(1,b,c)$ is semi-Gorenstein-projective if and only if  $\omega^s(1,b,c)\notin E$
for all $s\ge 0.$ Since $\omega^s(1,b,c) = (1,q^sb,c_s)$ for some $c_s\in k$,
we see that $M(1,b,c)$ is semi-Gorenstein-projective if and only if  $1+q^s\neq 0$
for all $s\ge 0$, thus if and only if  $q^{-s} \neq -b$ for all $s\ge 0$. Write $t = -s.$

Similarly, $M(1,b,c)$ is $\infty$-torsionfree if and only if  
$\omega^{-s}(1,b,c) \notin T$ for all $s\ge 0,$ thus if and only if  $1+q^{-1}q^{-s}b \neq 0$
for all $s\ge 0$, if and only if  $-b \neq q^{s+1}$ for all $s\ge 0.$ Write $t = s+1$. 
$\s$
	\medskip
{\bf Corollary.} {\it The module
$M(1,b,c)$ is Gorenstein-projective if and only if  $b\neq -q^t$ for all $t\in \Bbb Z.$}
	\bigskip 
{\bf (7.4)} {\it 
Any module $M(1,0,c)$ with $c\in k$ is Gorenstein-projective with 
$\Omega$-period $1$ or $2$.}
	\medskip
Proof. According to (6.2), the modules $M(1,0,c)$ are extensionless
and torsionless. Since $\omega(1,0,c) = (1,0,-c)$, we see that $M(1,0,0)$
has $\Omega$-period 1, and $M(1,0,c)$ with $c\neq 0$ has $\Omega$-period $2$
in case the characteristic of $k$ is different from $2$, otherwise its
$\Omega$-period is also $1$.  $\s$
	\bigskip

{\bf (7.5) Proposition.} {\it If $o(q) = \infty$, then any module of the form $M(1,b,c)$
is semi-Gorenstein-projective or $\infty$-torsionfree} (whereas the
modules of the form $M(0,b,c)$ are never semi-Gorenstein-projective nor $\infty$-torsionfree).
	\medskip
Proof. 	The first assertion follows immediately from (7.3), the additional assertion
in the bracket is a consequence of (5.1), (6.4) and (6.5).
$\s$

	\medskip
{\bf (7.6) Proposition.} 
{\it If $M(1,b,c)$ belongs to an $\Omega\mho$-component of the
form $\Bbb A_n$, then $o(q) = n.$}
	\medskip
Proof. We consider an $\Omega\mho$-component of type $\Bbb A_n$, say containing 
a  module $M$ which is not torsionless. 
Since $M$ belongs to $T$, we have $M = M(1,-q,c)$
and the component consists of the modules $M,\ \Omega M,\ \dots,\ \Omega^{n-1}M$.
In particular, $\omega^{n-1}(1,-q,c)$ belongs to $E$. 
Now $\Omega^{n-1}M = M(\omega^{n-1}(1,-q,c)) = M(1,-q^n,c')$ for some $c'$. 
Since $\Omega^{n-1}M$ is not extensionless, $(1,-q^n, c')$ belongs to $E$,
thus $-q^n = -1$. This shows that $q^n = 1.$ Finally, for $1\le t < n$, we have
$q^t \neq 1,$ since otherwise $\omega^{t-1}(1,-q,c)$ would belong to $E$. $\s$
	\medskip
{\bf Corollary.} {\it If $o(q) = \infty$, then
all the $\Omega\mho$-components in $H$ are cycles or of type $\Bbb Z,$ 
or $-\Bbb N$, or $\Bbb N$. 
Thus, any module in $H$ is semi-Gorenstein-projective or
$\infty$-torsionfree.}
	\medskip
For $o(q) = \infty$, there are the following 
$\Omega\mho$-components of the form $-\Bbb N$ and $\Bbb N$:
$$
{\beginpicture
    \setcoordinatesystem units <2.8cm,1.2cm>
\put{\beginpicture
\multiput{$\bullet$} at  1.5 0  2.5 0  /
\put{$\sssize \blacksquare$} at 3.5 0 
\arra{1.65 0}{1.6 0}
\arra{2.65 0}{2.6 0}
\setdashes <1mm>
\arr{2.4 0}{1.6 0}
\arr{3.4 0}{2.6 0}
\plot 1 0  1.4 0 /
\multiput{$\dots$} at 0.75 0 /
\put{$M(1,-\!q^3,c_3)$} at 1.5 -.3
\put{$M(1,-\!q^2,c_2)$} at 2.5 -.3
\put{$M(1,-\!q,c_1)$} at 3.5 -.3
\endpicture} at 0 0
\put{\beginpicture
\multiput{$\bullet$} at  1 0  2 0 /
\put{$\ssize \blacklozenge$} at 0 0 
\arra{0.15 0}{0.1 0}
\arra{1.15 0}{1.1 0}
\arra{2.15 0}{2.1 0}
\setdashes <1mm>
\arr{0.9 0}{0.1 0}
\arr{1.9 0}{1.1 0}
\arr{2.5 0}{2.1 0}
\multiput{$\dots$} at 2.75 0 /
\put{$M(1,-\!1,d_0)$} at 0 -.3
\put{$M(1,-\!q^{-1},d_{1})$} at 1 -.3
\put{$M(1,-\!q^{-2},d_{2})$} at 2 -.3
\endpicture} at 0 -.8
\endpicture}
$$
with  arbitrary elements
$c_0,d_1 \in k$ and $c_{t+1} = -\frac1{1-q^t}c_t$ for $t \ge 1$, whereas
$d_{t+1} = -(1-q^{-t})d_t$ for $t\ge 0$.
Of course, $(1,-q,c_1)\in T$ and $(1,-1,d_0)\in E$, thus the module $M(1,-q,c_1)$ is pivotal
semi-Gorenstein-projective, whereas $M(1,-1,d_0)$ is pivotal
$\infty$-torsionfree. 
	\medskip
{\bf (7.7) The case that $q$ has finite multiplicative order.}
{\it Now let $o(q) = n < \infty$. Then the modules $M(1,-q^t,c)$ with $0 \le t < n$ and $c\in k$
belong to $\Omega\mho$-components of the form $\Bbb A_n$.}
These $\Omega\mho$-components look as follows:
$$
{\beginpicture
    \setcoordinatesystem units <2.8cm,1.2cm>
\multiput{$\bullet$} at  1 0  2.5 0   /
\put{$\ssize \blacklozenge$} at 0 0 
\put{$\sssize \blacksquare$} at 3.5 0 
\arra{0.15 0}{0.1 0}
\arra{1.15 0}{1.1 0}
\arra{2.65 0}{2.6 0}
\setdashes <1mm>
\arr{0.9 0}{0.1 0}
\arr{1.4 0}{1.1 0}
\arr{3.4 0}{2.6 0}
\plot 2.1 0  2.4 0 /
\multiput{$\dots$} at 1.75 0 /
\put{$M(1,-\!1,c_n)$} at -.1 -.3
\put{$M(1,-\!q^{n-1},c_{n-1})$} at .9 -.3
\put{$M(1,-\!q^2,c_2)$} at 2.6 -.3
\put{$M(1,-\!q,c_1)$} at 3.6 -.3
\endpicture}
$$
with an arbitrary element 
$c_1 \in k$ and $c_{t+1} = -\frac1{1-q^t}c_t$ for $1\le t < n$ (of course, 
$(1,-1,c_n)\in E$ and $(1,-q,c_1)\in T$).

Corollary (7.3) asserts that the
remaining modules $M(1,b,c)$ (those with $-b\notin q^{\Bbb Z}$)
are Gorenstein-projective. 
	\bigskip
{\bf (7.9) Proof of Theorem (1.5).}
	\smallskip
{\bf Torsionless modules:}
According to (2.7), an indecomposable torsionless module is isomorphic to a
left ideal. Of course, $k$ is torsionless. According to (2.2), 
a 2-dimensional indecomposable left ideal is isomorphic to 
$\Lambda(x-y)$ or $\Lambda z$.
According to (2.3), a 3-dimensional indecomposable torsionless module has to be
local, thus it is 
of the form $M(a,b,c)$, and (6.1) says that $a(a+q^{-1}b)\neq 0$ or else
$M(a,b,c)$ is equal to $M(0,1,0)$ or to $M(0,0,1)$.
	\smallskip
{\bf Extensionless modules:} We show:
{\it An indecomposable module $M$ of dimension at most $3$ 
with simple socle is not extensionless.}

Of course, $\Ext^1(k,\Lambda) \neq 0$, since otherwise we would have $\Ext^1(X,\Lambda) = 0$
for all modules $X$. 

Let $I$ be an indecomposable module of length 2. A projective cover of $I$ as an
$\overline\Lambda$-module provides an exact sequence 
 $0 \to k^2 \to \overline\Lambda \to I \to 0$. We apply 
$\Hom_{\overline\Lambda}(-,J)$, where $J = \rad\Lambda$. 
We obtain the exact sequence
$$
 0 \to \Hom_{\overline\Lambda}(I,J) \to 
 \Hom_{\overline\Lambda}(\overline\Lambda,J) 
 \to \Hom_{\overline\Lambda}(k^2,J) \to 
 \Ext_{\overline\Lambda}^1(I,J) \to 0.
$$
Now, $\dim\Hom_{\overline\Lambda}(I,J) \ge  
\dim\Hom_{\overline\Lambda}(k,J) = 2,$ 
$\dim\Hom_{\overline\Lambda}(\overline\Lambda,J) = \dim J = 5,$ and
finally 
$\dim\Hom_{\overline\Lambda}(k^2,J) = 4,$ thus 
$\dim\Ext_{\overline\Lambda}^1(I,J) \ge 1.$ This shows that there exists 
a non-split
exact sequence $\epsilon\:0 \to J @>u>> E \to I \to 0$ with some $\overline\Lambda$-module $E$. 
The inclusion map $\iota\:J \to \Lambda$ yields an induced exact sequence
$\epsilon'\:0 \to \Lambda \to E' \to I \to 0$. Assume that 
 $\epsilon'$ splits. Then  we 
obtain a map $v\:E \to \Lambda$ such that $vu = \iota$. Now $E$ is an $\overline\Lambda$-module,
thus of Loewy length at most $2$. Therefore $v\:E \to \Lambda$ maps into $\rad\Lambda = J$,
thus $v = \iota v'$ for some $v'\:E \to J$. But $\iota v'u = vu = \iota$ implies that 
$v'u$ is the identity map of $E$, thus $\epsilon$ splits, a contradiction. The exact sequence
$\epsilon'$ shows that $\Ext_{\Lambda}^1(I,\Lambda) \neq 0.$ Thus $I$ is not extensionless. 

A similar proof shows that $\Ext^1(V,\Lambda) \neq 0$ for any 3-dimensional module $V$ with
simple socle. Again, we use that $V$ is an $\overline\Lambda$-module (see (1.3) Proposition 1),
thus we start with an exact sequence $0 \to k^5 \to \overline\Lambda^2 \to V \to 0.$

This completes the proof that an indecomposable module $M$ of dimension at most $3$ 
with simple socle is not extensionless. The remaining indecomposable modules of dimension
at most 3 are the modules of the form $M(1,b,c)$. According to (5.1) $M(1,b,c)$
is extensionless if and only if  $b\neq -1.$

	\smallskip
{\bf Reflexive modules:} We recall from Part I that a module $M$ is reflexive if and only if  both $M$ and 
$\mho M$ are
torsionless. We show:
{\it A module $M$ with simple socle is not reflexive.}
Assume that $M$ has simple socle and is torsionless. Since $M$ has simple
socle, there is an embedding $M \to {}_\Lambda\Lambda,$ say with cokernel $Q$.
The elements $yx$ and $zx$
cannot both belong to $u(M)$, since the socle of $u(M)$ is simple. If $yx\notin u(M)$, then
$yxQ\neq 0$, otherwise $zxQ\neq 0$.
Let $f\:M \to {}_\Lambda\Lambda^t$ be a minimal left 
$\add(\Lambda)$-approximation; its cokernel
is $\mho M$.
There is $u'\:{}_\Lambda\Lambda^t \to \Lambda$ with $u'f = u.$ The map $u'$ has to be
surjective, since otherwise $u'$ would vanish on the socle of ${}_\Lambda\Lambda^t$.
This implies that the map $\mho M \to Q$ induced by $u'$ is also surjective.
Since $\mho M$ is indecomposable, non-projective and not annihilated by $\rad^2\Lambda$,
$\mho M$ cannot be torsionless.

Let us assume that $M$ is reflexive and $\dim M\le 3$. It follows that $M$ has to be a torsionless
module with $\dim M = 3.$ 
Since also $\mho M$ has to be torsionless, (6.5) shows that the cases $M(0,1,0)$ and $M(0,0,1)$ are
not possible, thus $M$ is of the form $M(1,b,c)$ with $b\neq -q.$ Using (6.2) and (6.1), 
we see that we also must have $b \neq -q^2.$ Conversely, the same references show that all the
modules $M(1,b,c)$ with $b\neq -q^i$ for $i=1,2$ are reflexive. 
	\smallskip
{\bf Semi-Gorenstein-projective and $\infty$-torsionfree modules.}
The semi-Gorenstein-projective modules are extensionless, the $\infty$-torsionfree
modules are reflexive. The previous considerations therefore show that we only have to consider
the modules of the form $M(1,b,c)$. (7.3) provides the conditions on $b$ so that
$M(1,b,c)$ is semi-Gorenstein-projective, $\infty$-torsionfree, or
Gorenstein-projective. 

If  $M(1,b,c)$ is pivotal semi-Gorenstein-projective, then $M(1,b,c)$ is not
torsionless, thus $b = -q.$ If $M(1,-q,c)$ is semi-Gorenstein-projective, then
$-q \neq -q^{-s}$ for all $s\ge 0$, thus $q^{s+1} \neq 1$ for all $s\ge 0$. This means
that $o(q) = \infty.$ Of course, there is also the converse: if $o(q) = \infty$,
then $M(1,-q,c)$ is pivotal semi-Gorenstein-projective.

A similar argument shows that $M(1,b,c)$ is pivotal $\infty$-torsionfree
if and only if  $o(q) = \infty$ and $b = -1.$
$\s$
	\medskip
{\bf Remark.} It seems worthwhile to note that {\it the set of modules $M(1,b,c)$ with
$b,c\in k$ is a union of $\Omega\mho$-components.} 

	\bigskip\bigskip
{\bf 8. Right modules.}
	\medskip
Recall that we write $U'(a,b,c)$ instead of $U(a,b,c)$, if we consider $U(a,b,c)$ 
as a right ideal and that $M'(a,b,c) = \Lambda_\Lambda/U'(a,b,c)$.
	\medskip 
{\bf (8.1) Proposition.} {\it Let $(a,b,c)\neq 0.$ Then
$$
\Omega M' (a,b,c) \simeq \left\{
\matrix M' (\omega'(a,b,c))        &&\text{if}  & a \neq 0,\qquad\quad\, &\qquad\qquad & (1) \cr
        M' (0,0,1)                 &&\text{if}  &\ a =0,\ bc\neq 0,  && (2) \cr
        y\Lambda \oplus   zx\Lambda     &&\text{if}  & a = 0,\ c = 0, && (3)  \cr      
        z\Lambda \oplus   yx\Lambda     &&\text{if}  & a = 0,\ b = 0. && (4)  
\endmatrix
\right.
$$}

Proof. We have $\Omega M' (a,b,c) = U'(a,b,c)_\Lambda$. According to (2.8),
$U'(a,b,c)_\Lambda = (a,b,c)\Lambda$ if $a\neq 0$ or $bc\neq 0$, and $U'(0,1,0) =
y\Lambda\oplus zx\Lambda$, $U'(0,0,1) = z\Lambda\oplus yx\Lambda.$

Consider the map $\pi\:\Lambda_\Lambda \to U'(a,b,c)$
defined by $\pi(1) = (a,b,c)$. We assume that $a\neq 0$ or $bc\neq 0$, thus $\pi$
is surjective. If $a\neq 0$, the formula (3.1) (3) asserts that $\omega'(a,b,c)$
is in the kernel of $\pi$, thus $\pi$ yields an epimorphism
$M'(\omega'(a,b,c)) = \Lambda_\Lambda/\omega'(a,b,c)\Lambda \to U'(a,b,c)$.
Since this is a map between 3-dimensional modules, it has to be an isomorphism.

If $a = 0$ and $bc\neq 0$, we use formula (3.1) (4) in order to get similarly
an isomorphism $M'(0,0,1) = \Lambda_\Lambda/(0,0,1)\Lambda \to U'(0,b,c)$.
$\s$
	\bigskip
{\bf (8.2)} {\it If a $3$-dimensional indecomposable right module $N$ is torsionless
and no embedding $N \to \Lambda_\Lambda$ is a left $\add(\Lambda_\Lambda)$-approximation,
then $\mho N$ has Loewy length $3$ and is not torsionless.}
	\medskip
Proof. Let $\phi\:N \to \Lambda_\Lambda^t$ be a minimal left 
$\add(\Lambda_\Lambda)$-approximation of $N$. Since $N$ is torsionless, 
$\phi$ is a monomorphism. By assumption, we must have $t\ge 2.$ 
It follows that the cokernel $\mho N$ of $\phi$ is an indecomposable
right $\Lambda$-module of length $6t-3$ with top of length $t$. 
But an indecomposable right $\Lambda$-module of Loewy length at most $2$
with top of length $t\ge 2$ is a right  $\overline\Lambda$-module
of length at most $4t-1$. Thus $6t-3 \le 4t-1$, therefore $2t \le 2$, thus
$t\le 1$, a contradiction. This shows that $\mho N$ has Loewy length equal to $3$.
Of course, $\mho N$ is not projective. Since 
an indecomposable non-projective torsionless right $\Lambda$-module has
Loewy length at most $2$, we see that $\mho N$ cannot be torsionless.
$\s$
	\bigskip 
{\bf (8.3) The right modules $M'(0,b,c)$.} 
 {\it 
The only right module of the form $M'(0,b,c)$ which is torsionless is $M'(0,0,1).$ 
The right module $\mho M'(0,0,1)$ has
Loewy length $3$ and thus it is not torsionless.
No right module of the form $M'(0,b,c)$ is extensionless.}
	\medskip
Proof. Let $N = M'(0,b,c)$.

(a) If $N$ is torsionless, then $b = 0$ (thus $(0\col b\col c) = (0\col 0 \col 1)$).
Namely, According to (2.9), $M'(0,b,c)$ arises as a right ideal and (8.1) shows that
this happens only for $b = 0$.

(b) {\it No embedding $M'(0,0,1) \to \Lambda_\Lambda$ is a 
left $\add(\Lambda_\Lambda)$-approximation.}
Proof. Let $\phi\:M'(0,0,1) \to \Lambda_\Lambda$ be an embedding. According to (2.9),
the image of $\phi$ is of the form $U'(0,b,c)$ with $bc \neq 0.$
Now $M'(0,0,1)$ has a factor module isomorphic to $(0,0,1)\Lambda$,
thus there is $f\:M'(0,0,1) \to \Lambda_\Lambda$ with image $(0,0,1)\Lambda$.
If $\phi$ is a left $\add(\Lambda_\Lambda)$-approximation, then 
there exists $f':\Lambda_\Lambda \to \Lambda_\Lambda$ with $f = f'\phi$. 
The homomorphism $f'$ is the left multiplication by some element 
$\lambda$ in $\Lambda$. If $\lambda$ belongs to $\rad\Lambda$, then 
the image of $f'\phi$ is contained in $\rad^2\Lambda = \soc\Lambda$.
If $\lambda$ is invertible, then the image of $f'\phi$ is 3-dimensional.
In both cases, we get a contradiction, since the image of $f$ is $(0,0,1)\Lambda$,
thus $2$-dimensional and not contained in $\soc\Lambda$.

(c) It follows from (8.2) that $\mho M'(0,0,1)$ has Loewy length 3 and is
not torsionless.

(d) A right module of the form $M'(0,b,c)$ is never extensionless: either 
$\Omega M'(0,b,c)$ is decomposable, or else $\Omega M'(0,b,c) = M'(0,0,1)$ and 
according to (b), no embedding $M'(0,0,1) \to \Lambda_\Lambda$
is a left $\add(\Lambda_\Lambda)$-approximation.
$\s$
	\medskip
{\bf Reformulation.} {\it The right 
modules $M'(0,1,c)$ are singletons in the $\Omega\mho$-quiver. 
The right
module $M'(0,0,1)$ belongs to an $\Omega\mho$-component of the form $\Bbb A_2$:}
$$
{\beginpicture
    \setcoordinatesystem units <2cm,1cm>
\put{\beginpicture
\put{$\sssize \blacksquare$} at 1 0 
\put{$\ssize \blacklozenge$} at 0 0 
\arra{0.15 0}{0.1 0}
\setdashes <1mm>
\arr{0.9 0}{0.1 0}
\put{$M'(0,0,1)$} at -.2 -.3
\put{$\mho M'(0,0,1)$} at 1.2 -.3
\endpicture} at 0 0
\endpicture}
$$
	\medskip

{\bf (8.4) The right modules $M'(1,b,c)$ with $c\neq 0$.} 
	\medskip
{\bf Proposition.} {\it  Let $c\neq 0.$ 
The right module $M'(1,b,c)$ is torsionless if and only if  $b\neq -1$, and 
then $\mho M'(1,b,c) = M'(\omega(1,b,c))$. 
Let $c'\neq 0.$ 
The right module $M'(1,b',c')$ is extensionless if and only if  $b'\neq -q$, and then
$\Omega M'(1,b',c') = M'(\omega'(1,b',c'))$.}
	\medskip
{\bf Remark.} If $b\neq -1$ and $c\neq 0$, then $\omega(1,b,c) = (1,b',c')$
with $b' \neq -q$ and some $c'\neq 0$. 
If $b'\neq -q$, then $\omega'(1,b',c') = (1,b,c)$ with $b \neq -1$ and some $c\neq 0.$
Thus, the proposition provides $\Omega\mho$-sequences
$$
 0 \to M'(1,b,c) \to \Lambda_\Lambda \to M'(1,b',c') \to 0
$$
with $b\neq -1$ and $b'\neq -q$ (and both $c,c'$ being non-zero). 
Any triple $(1,b,c)$ with $b\neq -1$ and $c\neq 0$ occurs on the left
and given $(1,b,c)$, then we have $(1,b',c') = \omega(1,b,c)$  on the right.
Any triple $(1,b',c')$ with $b'\neq -q$ and $c'\neq 0$ occurs on the right
and given $(1,b',c')$, then we have $(1,b,c) = \omega'(1,b',c')$ on the left.
	\bigskip
Proof of Proposition. We follow closely the proof of (5.1) and (6.1). We always assume that
$c\neq 0$. As in (5.2) one sees that $M'(1,b,c)$ is extensionless if and only if  
the image of every homomorphism $U'(1,b,c) \to \Lambda_\Lambda$ is contained in
$U'(1,b,c)$.
	\smallskip
(a) {\it The module $M'(1,-q,c)$ is not extensionless.} Proof. According to (8.1),
we have $U'(1,-q,c') \simeq \Omega M'(1,-q,c') \simeq M'(\omega'(1,-q,c')) = M'(1,-1,0)$
for all $c'\in k.$ Thus, there is a homomorphism
$U'(1,-q,0) \to \Lambda_\Lambda$ with image $U'(1,-q,0)$ and 
this image $U'(1,-q,0)$ is not contained in $U'(1,-q,c).$
	\smallskip
(b) {\it If $b \neq -q,$ then the module $M'(1,b,c)$ is extensionless.} For the
proof, we need three assertions (b1), (b2) (b3).
Note that (8.1) asserts that $U'(1,b,c) \simeq \Omega M'(1,b,c) \simeq M'(\omega'(1,b,c))
= M'(1,q^{-1}b,c')$, where $\omega'(1,b,c) = (1,q^{-1}b,c')$. 
	\smallskip
(b1) {\it The only right ideal isomorphic to $U'(1,b,c)$ is $U'(1,b,c)$ itself.}
Proof. Let $V$ be a right ideal of $\Lambda_\Lambda$ which is isomorphic to $U'(1,b,c)$,
say $V = U'(a'',b'',c'')$ for some triple $(a'',b'',c'')$.
By (8.1), we have $U(a'',b'',c'') \simeq \Omega M'(a'',b'',c'') = M'(a'',q^{-1}b'',d)$,
where $\omega'(a'',b'',c'') = (a'',q^{-1}b'',d)$ 
for some $d$. We must have $a''\neq 0$, since $M(a'',q^{-1}b'',d) \simeq U'(1,b,c)
\simeq M'(1,q^{-1}b,c')$. Thus, we may assume that $a'' = 1$ and then
$M'(1,q^{-1}b'',d) \simeq M'(1,q^{-1}b,c')$ implies that $(1,q^{-1}b'',d) = (1,q^{-1}b,c')$.
In particular, we have $b'' = b\neq -q$. The equality
$\omega'(1,b'',c'') = \omega'(1,b,c)$ yields $(1,b'',c'') = (1,b,c)$, see Proposition (3.2).
Therefore $V = U(1,b'',c'') = U(1,b,c)$.
	\smallskip
(b2) {\it The right ideal $z\Lambda$ is not a factor module of $U'(1,b,c)$.}
Proof. The right ideal $z\Lambda$ is annihilated by $x-y$ and $z$, thus Corollary (5.3)
asserts that the image $I$ of any homomorphism $M'(1,b',c') \to z\Lambda$
is a factor module of $\Lambda/((1,b,c)\Lambda+(x-y)\Lambda +z\Lambda).$
Now $(x+by+cz)\Lambda + (x-y)\Lambda + z\Lambda = \rad \Lambda$, since $b\neq -1$,
thus $I$ is simple or zero.
	\smallskip
(b3) {\it The right ideal $y\Lambda$ is not a factor module of $U'(1,b,c)$.}
Proof. The right ideal $y\Lambda$ is annihilated by $y$ and $z$, thus Corollary (5.3)
asserts that the image $I$ of any homomorphism $M'(1,b',c') \to y\Lambda$
is a factor module of $\Lambda/((1,b,c)\Lambda+y\Lambda +z\Lambda).$
Now $(x+by+cz)\Lambda + y\Lambda + z\Lambda = \rad \Lambda$, since $b\neq -1$,
thus $I$ is simple or zero.
	\smallskip

The assertions (b1), (b2) and (b3) show: if $\phi$ is any homomorphism 
$U'(1,b,c) \to \Lambda_\Lambda$ and its image $I$ is of dimension at least 2, then 
$I$ is contained in $U'(1,b,c)$. Of course, if $I$ is 1-dimensional, then $I$ is
contained in $\soc\Lambda_\Lambda$ and $\soc\Lambda_\Lambda \subseteq U'(1,b,c)$. Thus,
we have obtained a proof of (b). In addition, (8.1) asserts that
$\Omega M'(1,b,c) \simeq M'(\omega'(1,b,c))$.
	\smallskip
(c) {\it If $b\neq -1$, then $M'(1,b,c)$ is torsionless and $\mho M'(1,b,c) =
M'(\omega(1,b,c))$.} Proof. Let $\omega(1,b,c) = (1,b',c').$ Then $b' = qb \neq -q$,
and $\omega'(1,b',c') = 
\omega'\omega(1,b,c) = (1,b,c)$ by Proposition (3.2). According to
(8.1), we have $\Omega M'(1,b',c') \simeq M'(\omega'(1,b',c')) = M'(1,b,c)$. 
This shows that $M'(1,b,c)$ is torsionless.According to (b), the module $M'(\omega(1,b,c))$ is extensionless, thus $\mho M'(1,b,c) =
M'(1,b',c') = M'(\omega(1,b,c))$.
	\smallskip
(d) {\it The modules $M'(1,-1,c)$ are not torsionless.} Proof. Assume, for the contrary, that
$M'(1,-1,c)$ is torsionless, thus isomorphic to $U'(a',b',c')$ for some $(a',b',c')$. 
According to (8.1), we must have $a' \neq 0$, thus we can assume that $a' = 1,$ 
and $(1,-1,c) = \omega'(1,b',c') = (1,q^{-1}b',-(1+q^{-1} b')c')$. It follows that
$b' = -q$ and therefore $c = -(1+q^{-1} b')c' = 0,$ a contradiction.

This completes the proof of (8.4). $\s$
	\bigskip
{\bf Reformulation.}
{\it The neighborhood of $M'(1,b,c)$ with $c\neq 0$ in the $\Omega\mho$-quiver
looks as follows:
$$
{\beginpicture
    \setcoordinatesystem units <2cm,.65cm>
\put{\beginpicture
    \setcoordinatesystem units <2cm,1.1cm>
\multiput{$\bullet$} at  1 0  /
\multiput{$\circ$} at  0 0  2 0 /
\arra{0.85 0}{0.9 0}
\arra{1.85 0}{1.9 0}
\setdashes <1mm>
\arr{0.1 0}{0.9 0}
\arr{1.1 0}{1.9 0}
\put{$M'(1,b,c)$} at 1 -.3
\put{$M'(\omega(1,b,c))$} at -.2 -.3
\put{$M'(\omega'(1,b,c))$} at 2.2 -.3
\multiput{$\dots$} at -.35 0  2.35 0 /
\multiput{} at -1 0  3 0 /
\endpicture} at 0 0
\put{\beginpicture
    \setcoordinatesystem units <2cm,1.1cm>
\put{$\ssize \blacklozenge$} at 1 0 
\put{$\circ$} at 0 0
\arra{0.85 0}{0.9 0}
\setdashes <1mm>
\arr{0.1 0}{0.9 0}
\put{$M'(1,-q,c)$} at 1 -.3
\put{$M'(\omega(1,-q,c))$} at -.2 -.3
\multiput{$\dots$} at -.35 0  /
\multiput{} at -1 0  3 0 /
\endpicture} at -.4 -1.5
\put{\beginpicture
    \setcoordinatesystem units <2cm,1.1cm>
\put{$\sssize \blacksquare$} at 1 0 
\put{$\circ$} at  2 0
\arra{1.85 0}{1.9 0}
\setdashes <1mm>
\arr{1.1 0}{1.9 0}
\put{$M'(1,-1,c)$} at 1 -.3
\put{$M'(\omega'(1,-1,c))$} at 2.2 -.3
\multiput{$\dots$} at  2.35 0 /
\multiput{} at -1 0  3 0 /
\endpicture} at 0.4 -3

\put{$b\notin\{-1,-q\}$} at 3.1 0
\put{$b = -q \neq -1$} at 3.1 -1.5
\put{$b = -1\neq -q$} at 3.1 -3
\endpicture}
$$
and $M'(1,b,c)$ is a singleton in the $\Omega\mho$-quiver if 
$q = 1$ and  $b=-1$.}
	\medskip
Note that we want to use a fixed index set $\Bbb P^2$ both for the 
(left) modules $M(a\col b \col c)$ 
and the right modules $M'(a\col b\col c)$, Since we have drawn the dashed arrows in
the $\Omega\mho$-quiver of the left $\Lambda$-modules from right to left, we
now have drawn the dashed arrows in the $\Omega\mho$-quiver of the right
$\Lambda$-modules from left to right. 
	\medskip
As in section 7, we see that the $\Omega\mho$-components of the modules
$M'(1,b,c)$ with $c\neq 0$ are cycles, or of type $\Bbb Z, \Bbb N$ or $-\Bbb N$
in case $o(q) = \infty$, and cycles or of type $\Bbb Z$ or $\Bbb A_n$ in case
$o(q) = n < \infty.$ 
	\medskip
{\it For $o(q) = \infty$, the right modules $M'(1,-1,c)$ 
with $c\neq 0$ are pivotal semi-Gorenstein-projective, and the right modules
$M'(1,-q,c)$ with $c\neq 0$ are pivotal $\infty$-torsionfree.}

	\bigskip
{\bf (8.5) The right modules $M'(1,b,0)$.} 
	\medskip
The right modules $M'(1,b,0)$ have been considered already in Part I: these are just
the right ideals $m_\alpha\Lambda$, where $m_\alpha = x-\alpha y$. Namely, we have
$$
     M'(1,b,0) = (x+qby)\Lambda = m_{-qb}\Lambda
$$
for all $b\in k.$ (Proof: We have $M'(1,b,0) 
= \Lambda_\Lambda/U'(1,b,0) = \Lambda_\Lambda/(x+by)\Lambda \simeq
(x+qby)\Lambda$, where we use that $(x+qby)(x+by) = 0$ and that
both right ideals $(x+by)\Lambda$
and $(x+qby)\Lambda$ are 3-dimensional, see (2.8).) 
	\medskip
Let us recall the results presented in Part I using the present notation:
	\medskip 
If $b\notin -q^{\Bbb Z}$, then $M'(1,b,0)$ is Gorenstein-projective and its
$\Omega\mho$-component looks as follows:
$$
{\beginpicture
    \setcoordinatesystem units <2.5cm,.7cm>
\multiput{$\bullet$} at 0 0  1 0  2 0  3 0  4 0  /
\multiput{$\cdots$} at -.6 0  4.6 0 /
\put{$M'(1,q^{2}b,0)$} at 0 -.5 
\put{$M'(1,qb,0)$} at 1 -.5
\put{$M'(1,b,0)$} at 2 -.5
\put{$M'(1,q^{-1}b,0)$} at 3 -.5
\put{$M'(1,q^{-2}b,0)$} at 4 -.5
\arra{-.15 0}{-.1 0}
\arra{0.85 0}{0.9 0}
\arra{1.85 0}{1.9 0}
\arra{2.85 0}{2.9 0}
\arra{3.85 0}{3.9 0}
\setdashes <1mm>
\arr{-.4 0}{-.1 0}
\arr{0.1 0}{0.9 0}
\arr{1.1 0}{1.9 0}
\arr{2.1 0}{2.9 0}
\arr{3.1 0}{3.9 0}
\plot 4.1 0 4.4 0 /
\endpicture}
$$
In particular, if $o(q) = n$, then these $\Omega\mho$-components are cycles with $n$ vertices, whereas for $o(q) = \infty$, one obtains $\Omega\mho$-components of type
$\Bbb Z$.

For $o(q) = \infty$, there are three remaining $\Omega\mho$-components:
$$
{\beginpicture
    \setcoordinatesystem units <2.5cm,.7cm>
\multiput{$\bullet$} at 0 0  3 0   4 0  /
\multiput{$\sssize\blacksquare$} at  1 1  2 1 /
\multiput{$\ssize\blacklozenge$} at 2 0  1 0  /
\multiput{$\cdots$} at -.3 0  4.3 0 /
\put{$M'(1,-q^{2},0)$} at 0 -.5 
\put{$M'(1,-q,0)$} at 1 -.5
\put{$M'(1,-1,0)$} at 2 -.5
\put{$M'(1,-q^{-1},0)$} at 3 -.5
\put{$M'(1,-q^{-2},0)$} at 4 -.5
\multiput{$\cdots$} at -.6 0  4.6 0 /
\arra{-.15 0}{-.1 0}
\arra{0.85 0}{0.9 0}
\arra{1.85 0.15}{1.9 0.1}
\arra{2.85 0.15}{2.9 0.1}
\arra{3.85 0}{3.9 0}

\setdashes <1mm>
\arr{-.4 0}{-.1 0}
\arr{0.1 0}{0.9 0}
\arr{1.1 0.9}{1.9 0.1}
\arr{2.1 0.9}{2.9 0.1}
\arr{3.1 0}{3.9 0}
\plot 4.1 0 4.4 0 /
\put{$\mho M'(1,-q^{-1},0)$} at 2.1 1.5
\put{$\mho M'(1,-1,0)$} at .9 1.5
\endpicture}
$$
These $\Omega\mho$-components are of type $\Bbb N, \Bbb A_2$ and $-\Bbb N$, respectively.

For $2\le n = o(q)< \infty$, there are two remaining $\Omega\mho$-components,
one is of type $\Bbb A_2$, the other of type $\Bbb A_n$:
$$
{\beginpicture
    \setcoordinatesystem units <2.2cm,.7cm>
\multiput{$\bullet$} at 3 0  4 0  5.5 0  /
\multiput{$\sssize\blacksquare$} at   2 1  1 1  /
\multiput{$\ssize\blacklozenge$} at 2 0  6.5 0   /
\multiput{$\cdots$} at  4.75 0   /
\put{$M'(1,-1,0)$} at 1.9 -.5
\put{$M'(1,-q^{n-1},0)$} at 3 -.5
\put{$M'(1,-q^{n-2},0)$} at 4.2 -.5
\put{$M'(1,-q^{2},0)$} at 5.4 -.5
\put{$M'(1,-q,0)$} at 6.5 -.5
\arra{1.85 0.15}{1.9 0.1}
\arra{2.85 0.15}{2.9 0.1}
\arra{3.85 0}{3.9 0}
\arra{6.35 0}{6.4 0}
\arra{5.35 0}{5.4 0}
\setdashes <1mm>
\arr{5 0}{5.4 0}
\arr{1.1 0.9}{1.9 0.1}
\arr{2.1 0.9}{2.9 0.1}
\arr{3.1 0}{3.9 0}
\arr{5.6 0}{6.4 0}
\put{$\mho M'(1,-q^{-1},0)$} at 2.1 1.5
\put{$\mho M'(1,-1,0)$} at .9 1.5
\plot 4.1 0  4.4 0 /
\endpicture}
$$

In case $q = 1$, 
there is only one additional $\Omega\mho$-component (of type $\Bbb A_2$), namely
$$
{\beginpicture
    \setcoordinatesystem units <2.3cm,.7cm>
\multiput{$\sssize\blacksquare$} at   2 1   /
\multiput{$\ssize\blacklozenge$} at 3 0  /
\put{$M'(1,-1,0)$} at 3 -.5
\arra{2.85 0.15}{2.9 0.1}
\setdashes <1mm>
\arr{2.1 0.9}{2.9 0.1}
\put{$\mho M'(1,-1,0)$} at 1.9 1.5
\endpicture}
$$
	\medskip
{\bf (8.6)} Similar to Theorem (1.5), here is the summary which characterizes the right modules of dimension at most 3 with relevant properties.
	\medskip
{\bf  Theorem.} {\it An indecomposable right module $N$ of dimension at most $3$ is 
\item{$\bullet$} torsionless if and only if  $N$ is simple or isomorphic to $y\Lambda$, to $z\Lambda$,
to a module $M'(1,b,c)$ with $b\neq -1$, to $M'(1,-1,0)$ or to $M'(0,0,1).$
\item{$\bullet$} extensionless if and only if  $N$ is isomorphic to a module $M'(1,b,c)$
   with $b\neq -q$;
\item{$\bullet$} reflexive if and only if  $M$ is isomorphic to a module $M'(1,b,c)$ with
   $b\neq -q^i$ for $i = -1, 0$;
\item{$\bullet$} Gorenstein-projective if and only if  $N$ is isomorphic to a module 
   $M'(1,b,c)$ with $b\neq -q^i$ for $i\in \Bbb Z$;
\item{$\bullet$} semi-Gorenstein-projective if and only if  $N$ is isomorphic to 
   a module $M'(1,b,c)$ with $b\neq -q^i$ for $i\ge 0$ or to 
   a module $M'(1,-1,c)$ with $c\neq 0$;
\item{$\bullet$} $\infty$-torsionfree if and only if  $N$ is isomorphic to 
    a module $M'(1,b,c)$ with $b\neq -q^i$ for $i\le 0$;
\item{$\bullet$} pivotal semi-Gorenstein-projective if and only if  $o(q) = \infty$ and 
   $N$ is isomorphic to a module $M'(1,-1,c)$ with $c\neq 0$;
\item{$\bullet$} pivotal $\infty$-torsionfree if and only if  $o(q) = \infty$ and 
   $N$ is isomorphic to a module $M'(1,-q,c)$.
\par}
$\s$
	\medskip
Whereas the set of modules $M(1,b,c)$ with 
$b,c\in k$ is a union of $\Omega\mho$-components, the right modules behave differently:
as we have seen already in Part I, 7.2, 
the $\Omega\mho$-component containing the right module $M(1,-1,0)$ consists of
$M(1,-1,0)$ and the 9-dimensional right module $\mho M(1,-1,0)$. 

	\bigskip\bigskip
{\bf 9. The $\Lambda$-dual of $M(1,b,c)$ and 
$M'(1,b,c)$.} 
	\medskip
We need the following (of course well-known) Lemma.
	\medskip 
{\bf (9.1) Lemma.} {\it Let $R$ be a ring and $w\in R$. If any 
left-module homomorphism $Rw \to {}_RR$ 
maps $w$ into $wR$, then $\Hom(Rw,{}_RR) \simeq wR$ as right $R$-modules.}
	\medskip
Proof. Let $u\:Rw\to {}_RR$ be the inclusion map. 
We have $\Hom(Rw,{}_RR) = uR,$ since for any homomorphism $f\:Rw \to {}_RR$, there
is $\lambda\in R$ with  $f(w) = w\lambda$, thus  $f = u\lambda.$
Now $I = \{r\in R\mid wr = 0\}$ is a right ideal and $R_R/I \simeq wR$ 
as right modules (an isomorphism is given by the map 
$R_R \to wR$ defined by $1\mapsto w$). Since $I = \{r\in R\mid
ur = 0\}$, we have in the same way $R_R/I \simeq uR,$ and therefore 
$wR \simeq R_R/I \simeq uR = \Hom(Rw,{}_RR).$
$\s$

	\medskip
{\bf (9.2) Lemma.} {\it If $(1,b,c)$ is different from $(1,-1,0)$, then $M'(1,b,c) \simeq
\Tr M(1,b,c)$ and $M(1,b,c) \simeq \Tr M'(1,b,c)$.}
	\medskip
Proof. We have $U'(1,b,c) = (1,b,c)\Lambda$, and since $(1,b,c) \neq (1,-1,0)$, we also
have $U(1,b,c) = \Lambda(1,b,c)$.
By definition, $M(1,b,c) = {}_\Lambda\Lambda/U(1,b,c)$, thus
$M(1,b,c)$ is the cokernel of the right multiplication
$r_{(1,b,c)}\:{}_\Lambda\Lambda \to {}_\Lambda\Lambda$ and 
$\Tr M(1,b,c)$ is the cokernel of the left multiplication 
$l_{(1,b,c)}\:\Lambda_\Lambda \to
\Lambda_\Lambda$,
thus isomorphic to $\Lambda_\Lambda/(1,b,c)\Lambda = \Lambda_\Lambda/U'(1,b,c)$. 
$\s$
	\medskip
{\bf (9.3) Proposition.}
{\it If $b\notin\{-q,-q^2\}$, then $M(1,b,c)$ is reflexive and 
$$
 M(1,b,c)^* = M'((\omega')^2(1,b,c)).
$$
If $b\notin\{-1,-q^{-1}\}$, then $M'(1,b,c)$ is reflexive and }
$$
 M'(1,b,c)^* = M(\omega^2(1,b,c)).
$$
	\medskip
Proof. According to (7.1), we have the following two $\Omega\mho$-sequences:
$$
\gather
 0 \to M(1,b,c) \to {}_\Lambda\Lambda \to M(\omega'(1,b,c)) \to 0,\cr
 0 \to M(\omega'(1,b,c)) \to {}_\Lambda\Lambda \to M((\omega')^2(1,b,c)) \to 0
\endgather
$$
(the first one, since $\omega'(1,b,c) = (1,b',c')$ with $b' = q^{-1}b \neq -1;$
the second one, since $(\omega')^2(1,b,c) = (1,b'',c'')$ with $b'' = q^{-2}b \neq -1$)
This implies that $M(1,b,c)$ is reflexive and that
$X = \mho^2 M(1,b,c) = M((\omega')^2(1,b,c))$ is a module with $\Ext^i(X,\Lambda) = 0$ for
$i=1,2.$ According to Part I, Lemma 2.5, we have $\Tr X = (\Omega^2 X)^*.$
On the one hand, $\Omega^2 X = \Omega\mho M(1,b,c) = M(1,b,c)$.
On the other hand, (9.2) shows that 
$\Tr X = \Tr M((\omega')^2(1,b,c)) = M'((\omega')^2(1,b,c))$, since
$(\omega') ^2(1,b,c) = (1,q^{-2}b,c'')$ for some $c''$ and $q^{-2}b \neq -1$.
This yields the first assertion. 
The second can be shown in the same way, or just by applying 
the $\Lambda$-duality to $M(1,b,c)^* = M'((\omega')^2(1,b,c))$.
$\s$
	\medskip
{\bf (9.4) Proposition.} {\it For all $b,c\in k$,}
$$
 M(1,b,c)^* = M'((\omega')^2(1,b,c)).
$$
	\medskip
In particular, for all $b,c\in k$, the right module $M(1,b,c)^*$ is again 3-dimensional
and local. 
	\medskip 
Whereas $(\omega')^2$ is a bijection from $\{(1,b,c)\mid b\notin\{-q,-q^2\}\}$ onto
$\{(1,b,c)\mid b\notin\{-1,-q^{-1}\}\}$, we should stress that $(\omega')^2(1,-q,c) = 
(1,-q^{-1},0)$ and that $(\omega')^2(1,-q^2,c) = (1,-1,0)$ for all $c\in k.$ Thus, 
(9.3) combines the first assertion of (9.2) with the corresponding assertion
for the remaining cases, namely:
$$
 M(1,-q,c)^* = M'(1,-q^{-1},0) \quad{and}\quad M(1,-q^2,c)^* = M'(1,-1,0),
$$
for all $c\in k$.

	\medskip
Proof of Proposition. According to (9.2), we only have to consider the cases where
$b = -q$ or $b= -q^2.$ 
	\smallskip
Case 1. Let $b = -q.$ 
As we have seen in (6.2), the module $M(1,-q,c)$ is not torsionless.
Now obviously, there is a surjective homomorphism $M(1,-q,c) \to \Lambda(1,-1,0)$
with kernel $zM(1,-q,c)$. It follows that $zM(1,-q,c)$ is contained in the kernel
of every homomorphism $M(1,-q,c) \to {}_\Lambda\Lambda$ and therefore 
$M(1,-q,c)^* = (\Lambda(1,-1,0))^*$. Now, 
$(\Lambda(1,-1,0))^* \simeq (1,-1,0)\Lambda = U'(1,-1,0)$, 
as shown in Part I, 6.5.
On the other hand, according to (8.1), we have
$U'(1,-1,0) = \Omega M'(1,-1,0) = M'(\omega'(1,-1,0))$ and 
$\omega'(1,-1,0) = (1,-q^{-1},0).$
	\smallskip
Case 2: $b = -q^2$ and $o(q) = 2.$
The assumption $o(q) = 2$ means that $q = -1\neq 1$, in particular,
the characteristic of $k$ is different from $2$, and we have $b = -1.$ 
Since $q = -1$ and the characteristic of $k$ is different from $2$, 
(4.1) asserts that 
$$
 \Lambda(1,1,-2c) = U(1,1,-2c) = \Omega M(1,1,-2c) = M(\omega(1,1,-2c)) = M(1,-1,c).
$$
On the other hand, we have
$$
 (1,1,-2c)\Lambda = U'(1,1,-2c) = \Omega M'(1,1,-2c) = M'(\omega'(1,1,-2c)) = 
  M'(1,-1,0).
$$
We claim that any homomorphism $\Lambda(1,1,-2c) \to {}_\Lambda\Lambda$ maps
$(1,1,-2c)$ into $(1,1,-2c)\Lambda$. Namely, let $\phi\:\Lambda(1,1,-2c) \to {}_\Lambda\Lambda$ be a homomorphism. 
Now $\Lambda(1,1,-2c)$ is 3-dimensional, thus equal to $U(1,1,-2c)$,
and ${}_\Lambda\Lambda/U(1,1,-2c) \simeq M(1,1,-2c).$ According to (5.1), the module
$M(1,1,-2c)$ is extensionless, since $1+1\neq 0$. The implication (i) to (iv)
in (5.2) shows that $\phi(1,1,-2c) \in (1,1,-2c)\Lambda$. 

Since any homomorphism $\Lambda(1,1,-2c) \to {}_\Lambda\Lambda$ maps
$(1,1,-2c)$ into $(1,1,-2c)\Lambda$, Lemma (9.0) implies that the right modules
$(\Lambda(1,1,-2c))^*$ and $(1,1,-2c)\Lambda$ are isomorphic, thus 
$M(1,-1,c)^* \simeq M'(1,-1,0).$
	\medskip
Case 3. $b = -q^2$ and $o(q) \ge 3$. 
There is the $\Omega\mho$-sequence
$$
 \epsilon\:\qquad  0 \to M(1,-q^3,c') \to {}_\Lambda\Lambda \to M(1,-q^2,c) \to 0
$$
for some $c'$ (here we use that $q^2\neq 1$). 
The $\Lambda$-dual of $\epsilon$ is the exact sequence
$$
 0 \to M(1,-q^2,c)^* \to \Lambda_\Lambda \to M(1,-q^3,c')^* \to 0.
$$
Since $q^2 \neq 1$,  proposition (9.3) asserts that 
$M(1,-q^3,c')^* = M'(1,-q,c'')$ for some $c''$. Altogether we see that
$$
 M(1,-q^2,c)^* \simeq \Omega (M(1,-q^3,c')^*) = \Omega M'(1,-q,c'') \simeq M'(1,-1,0),
$$
where the final isomorphism is due to (8.1).
$\s$
	\medskip
{\bf (9.5)}  The algebra $\Lambda = 
\Lambda(q)$ with $o(q) = \infty$ was exhibited in Part I in order to present 
a module $M$  which is
not torsionless, such that $M$ and $M^*$ both are semi-Gorenstein-projective:
namely the module $M = M(1,-q,0)$ with $M^* = M'(1,-q,0)$.
Now we see: {\it all the modules $M(1,-q,c)$ with $c\in k$ are modules which are
semi-Gorenstein-projective and not torsionless, and that the 
$\Lambda$-duals $M(1,-q,c)^* \simeq M'(1,-q^{-1},0)$ are semi-Gorenstein-projective.}
We should stress that this concerns a 1-parameter family $M(1,-q,c)$ (with $c\in k$)
of semi-Gorenstein-projective left modules, and the single semi-Gorenstein-projective
right module $M(1,-q^{-1},0)$.
	\medskip
{\bf (9.6) Proposition.} {\it Let $b,c\in k.$}
$$
 M'(1,b,c)^* =
  \left\{
 \matrix M(\omega^2(1.b.c)) & \text{if} &\ \, b\notin\{-1,-q^{-1}\},\qquad \qquad  \cr
  U(0,0,1) & \text{if}                  & b = -1,\ c \neq 0, \qquad \qquad \cr
  U(1,-q,0)+U(0,0,1) & \text{if}        & b = -1,\ c = 0,\qquad \qquad \cr
  M(0,0,1) & \text{if}                  & \, b = -q^{-1},\ c \neq 0,\ q\neq 1,\cr
  U(1,-1,0) & \text{if}                 & \, b = -q^{-1},\ c = 0,\ q\neq 1.
 \endmatrix \right.
$$
	\medskip
Whereas we saw in (9.4) that all the right modules $M(1,b,c)^*$ are 3-dimensional and local,
not all the modules $M'(1,b,c)^*$ are 3-dimensional and local: the module 
$M'(1,-1,0)^* =
U(1,-q,0)+U(0,0,1)$ has dimension 4, whereas the modules 
$M'(1,-1,c)^* = U(0,0,1)$ for $c\neq 0$
and, in case $q\neq 1$, the module $M'(1,-q^{-1},0)^* = U(1,-1,0)$ are decomposable. 
	\medskip
Proof. According to (9.3), we only have to deal with the cases with $b \in\{-1,-q^{-1}\}.$
If $c = 0$, then we can refer to Part I. 
For $b = -1$, the end of 7.1 in Part I shows that $M'(1,-1,0)^* \simeq M(1,-q^2,0)^{**} \simeq 
U(1,-q,0)+U(0,0,1)$. For $b = -q^{-1} \neq -1$, the end of 6.7 in Part I asserts that
$M'(1,-q^{-1},0)^* \simeq (M(1,-q,0)^{**} \simeq \Omega M(1,-1,0) \simeq U(1,-1,0).$
	\smallskip 
Now, we assume that $c\neq 0$. As in the proof of (9.4), we consider again 3 cases.
	\smallskip
Case 1. $b = -1$. 
The module $M'(1,-1,c)$ with $c\neq 0$ is not torsionless, see (8.4). Since 
the factor module  $M'(1,-1,c)/M'(1,-1,c)z$ is isomorphic to $(0,0,1)\Lambda$,
it follows that $M'(1,-1,c)^* \simeq ((0,0,1)\Lambda)^*$ and an easy calculation yields 
$((0,0,1)\Lambda)^* \simeq U(0,0,1)$. Namely, 
the inclusion map $u\:z\Lambda \to \Lambda_\Lambda$ satisfies $yu = 0$ and $zu = 0$,
thus a basis of $(z\Lambda)^*$ is given by $u,\ xu$ and the map $f\:z\Lambda
\to \Lambda_\Lambda$ with $f(z) = yx$, so that $(z\Lambda)^* \simeq 
{}_\Lambda\Lambda/(\Lambda y+\Lambda z) \oplus k \simeq U(0,0,1).$
	\smallskip
Case 2. $b = -q^{-1}$ and $o(q) = 2.$ 
Thus, the characteristic of $k$ is different from $2$, $q = -1$ and $b = 1$.
The module $M'(1,1,c)$ is torsionless: namely, by (8.1) we have
$M'(1,1,c) \simeq \Omega M'(1,-1,-\frac c2),$ since $\omega'(1,-1,-\frac c2) = (1,1,c)$.
Now, $\Omega M'(1,-1,-\frac c2) \simeq U'(1,-1,\frac c2) = (1,-1,\frac c2)\Lambda$.
Since $q\neq 1$, the right module 
$M'(1,-1,-\frac c2)$ is extensionless by (8.4), thus we can use
(5.2) and (9.1) in order to see that $((1,-1,\frac c2)\Lambda)^* \simeq \Lambda(1,-1,\frac c2).$
By (4.1) (2), we have $\Lambda(1,-1,\frac c2) = U(1,-1,\frac c2) \simeq 
\Omega M((1,-1,-\frac c2)) \simeq M(0,0,1).$
	\smallskip 
Case 3. $b = -q^{-1}$ and $o(q) \ge 3.$ 
There is the $\Omega\mho$-sequence
$$
 0 \to M'(1,-q^{-2},c') \to \Lambda_\Lambda \to M'(1,-q^{-1},c) \to 0
$$
for $c' = \lambda c$ with $\lambda\neq 0$
(here we use that $q^2\neq 1$). The $\Lambda$-dual is the exact sequence
$$
 0 \to M'(1,-q^{-1},c)^* \to {}_\Lambda\Lambda \to M'(1,-q^{-2},c')^* \to 0.
$$
We assume that $q\neq 1$ and $q\neq 2$. Then by Proposition (9.2), we have
$M'(1,-q^{-2},c')^* = M(1,-1,c'')$ for some multiple $c'' = \lambda'c'$ with
$\lambda'\neq 0$. It follows that
$M'(1,-q^{-1},c)^* = \Omega M(1,-1,c'')$ and $c'' = 0$ if and only if  $c = 0.$ By (4.1),
we have $\Omega M(1,-1,c'') = M(0,0,1)$ in case $c\neq 0$, and
$\Omega M(1,-1,0) = U(1,-1,0)$  in case $c = 0.$ $\s$
	\medskip
{\bf (9.7) Corollary.} {\it Let $N$ be a right $\Lambda$-module of dimension at
most $3$ which is semi-Gorenstein-projective, but not Gorenstein-projective. 
Then $N^*$ is not semi-Gorenstein-projective.}
	\smallskip
Proof. 
According to (8.6), $N$ is isomorphic to a right module of the form $M'(1,-q^i,c)$
with $i\le -1$ and $c\in k$ or of the form $M'(1,-1,c)$ with $c\neq 0.$ We
apply (9.6).
If $i\le -2$, then $N^* = M'(1,-q^i,c)^* = M(1,-q^{i+2},c')$ for some $c'$,
and according to (1.5), $N^*$ is not semi-Gorenstein-projective, since $i+2 \le 0$.
If $i = -1$, then $N^*$ is isomorphic to $M(0,0,1)$ or to $U(1,-1,0)$.
If $N = M'(1,-1,c)$ with $c\neq 0,$ then $N^*$ is isomorphic to $U(0,0,1).$ 
But by (1.5), $M(0,0,1),$ $U(1,-1,0)$ and $U(0,0,1)$ are not semi-Gorenstein-projective.
$\s$
	\bigskip
{\bf 10. The general context.}
	\medskip
Our detailed study of the algebra $\Lambda(q)$ in Part I and Part II should
be seen in the frame of looking at Gorenstein-projective (or, more general,
semi-Gorenstein-projective and $\infty$-torsionfree modules) over short local
algebras. 
	
Let $A$ be a finite-dimensional local $k$-algebra with radical $J$ such that
$A/J = k$. Such an algebra is said to be {\it short} provided $J^3 = 0.$ 
In commutative ring theory, the short local algebras have attracted a lot of
interest, since some conjectures have been disproved by  
looking at modules over short algebras, see [AI\c S] for a corresponding account.

Let us assume now that $A$ is short, but not necessarily commutative.
Let $e = \dim J/J^2$ and $a = \dim J^2$
(thus $0 \le a \le e^2$). If there exists an indecomposable module which is
semi-Gorenstein-projective or $\infty$-torsionfree, 
but not projective, then either $A$ is self-injective,
so that $a\le 1$ (and $e =1$ in case $a = 0$), or else $a = e-1$ and 
$J^2 = \soc {}_AA = \soc A_A$. 

Of course, if $A$ is self-injective, then all
modules are Gorenstein-projective, thus the interesting case is the case $a = e-1$.
Our algebra $\Lambda(q)$ is of this kind (with $a = 2$), as is the Jorgensen-\c Sega
algebra [J\c S] (with $a = 3$). 

Not only the shape of the algebras is very
restricted, also the modules themselves are very special: Let $A$ be a short local algebra
which is not self-injective. Let $M$ be indecomposable and not projective.
If $M$ is semi-Gorenstein-projective and torsionless, or if $M$ is $\infty$-torsionfree,
(in particular, if $M$ is Gorenstein-projective), 
then $\soc M = \rad M$ and $\dim \soc M = a\cdot \dim \top M$ (by definition, 
$\top M = M/\soc M$). Also, if $M$ is semi-Gorenstein-projective and torsionless, then 
$\dim \Omega^i M = \dim M$ for all $i\in \Bbb N$, whereas 
if $M$ is $\infty$-torsionfree, then $\dim \mho^i M = \dim M$ for all $i\in \Bbb N$.
These assertions have been shown by Christensen
and Veliche in the case that $A$ is commutative, see [CV], but actually the proofs 
do not have to be modified in the general case. There is an essential
difference between the commutative and the non-commutative algebras: If $A$ is
commutative, then all local modules which are semi-Gorenstein-projective or
$\infty$-torsionfree are Gorenstein-projective, whereas this is not true for $A$
non-commutative. For a general discussion, we refer to [RZ2] (and we have to thank D.
Jorgensen for his advice concerning the present knowledge in the commutative case).

Thus, for our algebra $\Lambda(q)$, the non-projective
indecomposable modules which are semi-Gorenstein-projective and torsionless, 
or which are $\infty$-torsionfree, are
of dimension $3t$ with socle of dimension $2t$, where $t = \dim \top M$. 
For $t=1$, we deal with local modules with 2-dimensional socle: these are
precisely the modules studied in the present paper.
	\bigskip\bigskip

{\bf Appendix. A diagrammatic description of the modules $M(a\col b\col c)$.}
     \medskip
 If $M$ is a left $\Lambda$-module annihilated by $\rad^2\Lambda$, then
 it is a left $\overline\Lambda$-module. Since $\overline\Lambda$ is a commutative
 $k$-algebra, also $D(M) = \Hom(M,k)$ is a left $\overline\Lambda$-module, thus
 a left $\Lambda$-module.
      \medskip
 {\bf Proposition.} {\it Let $M$ be an indecomposable $3$-dimensional left
 $\Lambda$-module. Then $M$ or $D(M)$ is isomorphic to one of the
 following pairwise non-isomorphic $\overline\Lambda$-modules $M(a,b,c)$:}
 $$
 {\beginpicture
 \setcoordinatesystem units <1.9cm,2cm>

 \put{Case} at -1 -.3
 \put{(1)} at -1 -1
 \put{(2)} at -1 -2
 \put{(3)} at -1 -3
 \put{(4)} at -1 -4
 \put{(5)} at -1 -5
 \put{(6)} at -1 -6
 \put{(7)} at -1 -7

 \put{Position in $\Bbb P^2$} at 1.4 -.3
\put
 {\beginpicture
 \setcoordinatesystem units <.7cm,.7cm>
 \setdots <.5mm>
 \multiput{} at 0 0  2 1 /
 \plot 0 0  1 1  2 0  0 0 /
 \put{$\bullet$} at 1 1 
 \endpicture} at 1.4 -1
\put
 {\beginpicture
 \setcoordinatesystem units <.7cm,.7cm>
 \setdots <.5mm>
 \multiput{} at 0 0  2 1 /
 \plot 0 0  1 1  2 0  0 0 /
 \put{$\bullet$} at 2 0 
 \endpicture} at 1.4 -2
\put
 {\beginpicture
 \setcoordinatesystem units <.7cm,.7cm>
 \setdots <.5mm>
 \multiput{} at 0 0  2 1 /
 \plot 0 0  1 1  2 0  0 0 /
 \put{$\bullet$} at 0 0  
 \endpicture} at 1.4 -3
\put
 {\beginpicture
 \setcoordinatesystem units <.7cm,.7cm>
 \setdots <.5mm>
 \multiput{} at 0 0  2 1 /
 \plot 0 0  1 1  2 0  0 0 /
\setsolid
\plot 0.15 0  1.85 0 /
\plot 0.15 0.025  1.85 0.025 /
\plot 0.15 -.025  1.85 -.025 /
 \endpicture} at 1.4 -4
\put
 {\beginpicture
 \setcoordinatesystem units <.7cm,.7cm>
 \setdots <.5mm>
 \multiput{} at 0 0  2 1 /
 \plot 0 0  1 1  2 0  0 0 /
\setsolid
\plot 0.1 0.1  0.9 0.9 /
\plot 0.08 0.12  0.88 0.92 /
\plot 0.12 0.08  0.92 0.88 /
 \endpicture} at 1.4 -5
\put
 {\beginpicture
 \setcoordinatesystem units <.7cm,.7cm>
 \setdots <.5mm>
 \multiput{} at 0 0  2 1 /
 \plot 0 0  1 1  2 0  0 0 /
\setsolid
\plot 1.1 0.9  1.9 0.1 /
\plot 1.12 0.92  1.92 0.12 /
\plot 1.08 0.88  1.88 0.08 /
 \endpicture} at 1.4 -6
\put
 {\beginpicture
 \setcoordinatesystem units <.7cm,.7cm>
 \setdots <.4mm>
 \multiput{} at 0 0  2 1 /
 \plot 0 0  1 1  2 0  0 0 /
 \setshadegrid span <.2mm>
 \vshade 0.15 0.08 0.1  <,z,,> 1 0.08 .9  <z,,,> 1.85 0.08 0.1 /

 \endpicture} at 1.4 -7

\put{Modules} at 0 -.3
\put{$M(0, 0 , 1)$} at 0 -1
\put{$M(0, 1 , 0)$} at 0 -2
\put{$M(1, 0 , 0)$} at 0 -3

\put{$M(1, b , 0)$} at 0 -4
\put{$b\in k^*$} at 0 -4.2
\put{$M(1, 0 , c)$} at 0 -5
\put{$c\in k^*$} at 0 -5.2
\put{$M(0, 1 , c)$} at 0 -6
\put{$c\in k^*$} at 0 -6.2

\put{$M(1, b , c)$} at 0 -7
\put{$b, c\in k^*$} at 0 -7.2

\put{Diagram} at 3 -.3

\put{\beginpicture
 \setcoordinatesystem units <1cm,1cm>
 \put{$v$} at 1 1
 \put{$v_1$} at 0 0
 \put{$v_2$} at 2 0
 \arr{0.8 0.8}{0.2 0.2}
 \arr{1.2 0.8}{1.8 0.2}
 \put{$\ssize x$} at 0.4 0.6
 \put{$\ssize y$} at 1.6 0.6
 \put{} at 1 -0.3
 \endpicture} at 3 -1
\put{\beginpicture
 \setcoordinatesystem units <1cm,1cm>
 \put{$v$} at 1 1
 \put{$v_1$} at 0 0
 \put{$v_2$} at 2 0
 \arr{0.8 0.8}{0.2 0.2}
 \arr{1.2 0.8}{1.8 0.2}
 \put{$\ssize x$} at 0.4 0.6
 \put{$\ssize z$} at 1.6 0.6
 \put{} at 1 -0.3
 \endpicture} at 3 -2
\put{\beginpicture
 \setcoordinatesystem units <1cm,1cm>
 \put{$v$} at 1 1
 \put{$v_1$} at 0 0
 \put{$v_2$} at 2 0
 \arr{0.8 0.8}{0.2 0.2}
 \arr{1.2 0.8}{1.8 0.2}
 \put{$\ssize y$} at 0.4 0.6
 \put{$\ssize z$} at 1.6 0.6
 \put{} at 1 -0.3
 \endpicture} at 3 -3

\put{\beginpicture
 \setcoordinatesystem units <1cm,1cm>
 \put{$v$} at 1 1
 \put{$v_1$} at 0 0
 \put{$v_2$} at 2 0
 \arr{0.8 0.8}{0.2 0.2}
 \arr{1.2 0.8}{1.8 0.2}
 \put{$\ssize y$} at 0.4 0.6
 \put{$\ssize z$} at  1.6 0.6
\setdashes <1mm>
\setquadratic
\plot 0.93 0.83  0.7 0.3  0.2 0.1 /
\arr{0.25 0.11}{0.2 0.1}
\put{$\ssize x$} at 0.8 0.15
 \put{with $xv = -bv_1$} at 1 -0.35
 \endpicture} at 3 -4

\put{\beginpicture
 \setcoordinatesystem units <1cm,1cm>
 \put{$v$} at 1 1
 \put{$v_1$} at 0 0
 \put{$v_2$} at 2 0
 \arr{0.8 0.8}{0.2 0.2}
 \arr{1.2 0.8}{1.8 0.2}
 \put{$\ssize y$} at 0.4 0.6
 \put{$\ssize z$} at  1.6 0.6
\setdashes <1mm>
\setquadratic
\plot 1.07 0.83  1.3 0.3  1.8 0.1 /
\arr{1.75 0.11}{1.8 0.1}
\put{$\ssize x$} at 1.2 0.15
 \put{with $xv = -cv_2$} at 1 -0.35
 \endpicture} at 3 -5

\put{\beginpicture
 \setcoordinatesystem units <1cm,1cm>
 \put{$v$} at 1 1
 \put{$v_1$} at 0 0
 \put{$v_2$} at 2 0
 \arr{0.8 0.8}{0.2 0.2}
 \arr{1.2 0.8}{1.8 0.2}
 \put{$\ssize x$} at 0.4 0.6
 \put{$\ssize z$} at 1.6 0.6
\setdashes <1mm>
\setquadratic
\plot 1.07 0.83  1.3 0.3  1.8 0.1 /
\arr{1.75 0.11}{1.8 0.1}
\put{$\ssize y$} at 1.2 0.15
 \put{with $yv = -cv_2$} at 1 -0.35
 \endpicture} at 3 -6

\put{\beginpicture
 \setcoordinatesystem units <1cm,1cm>
 \put{$v$} at 1 1
 \put{$v_1$} at 0 0
 \put{$v_2$} at 2 0
 \arr{0.8 0.8}{0.2 0.2}
 \arr{1.2 0.8}{1.8 0.2}
 \put{$\ssize y$} at 0.4 0.6
 \put{$\ssize z$} at  1.6 0.6
\setdashes <1mm>
\setquadratic
\plot 0.93 0.83  0.7 0.3  0.2 0.1 /
\arr{0.25 0.11}{0.2 0.1}
\put{$\ssize x$} at 0.8 0.15
\setdashes <1mm>
\plot 1.07 0.83  1.3 0.3  1.8 0.1 /
\arr{1.75 0.11}{1.8 0.1}
\put{$\ssize x$} at 1.2 0.15

 \put{with $xv = -bv_1-cv_2$} at 1 -0.35
 \endpicture} at 3 -7

\put{Characterization} at 5 -.3
\put{$zM = 0$} at 5 -1
\put{$yM = 0$} at 5 -2
\put{$xM = 0$} at 5 -3
\put{$xM = yM$} at 5 -4
\put{$xM = zM$} at 5 -5
\put{$yM = zM$} at 5 -6
\put{$xM,\,yM,\,zM$ non-zero} at 5 -6.9
\put{and pairwise different}  at 5 -7.1
\endpicture}
$$

The diagrams describe the modules $M = M(a,b,c)$ as follows:
The elements $v,v_1,v_2$ form a basis of $M$. Both elements $v_1,v_2$ are annihilated by
$x,y,z$. If there is drawn a solid arrow $v$
\beginpicture
    \setcoordinatesystem units <.5cm,.8cm>
\multiput{} at 0 0  1 0 /
\arr{0 0.1}{1 0.1}
\endpicture\
$v_i$ with $i\in\{1,2\}$ and with label $\alpha\in \{x,y,z\}$,
then $\alpha v = v_i.$ If there is a dashed arrow $v$
\beginpicture
    \setcoordinatesystem units <.5cm,.8cm>
\multiput{} at 0 0  1 0 /
\setdashes <1mm>
\arr{0 0.1}{1 0.1}
\endpicture\
$v_i$ with label $\alpha$, then $\alpha v = c_1v_1+c_2v_2$ with $c_i\neq 0$
(and we provide the coefficients $c_1,c_2$ below the diagram). Finally, $zv = 0$ in case (1),
$yv = 0$ in case (2), $xv = 0$ in case (3). 

The last column provides a characterization of the corresponding modules $M(a,b,c)$: 
For example, a local 3-dimensional $\Lambda$-module $M$ is a case-(1)-module 
provided $zM = 0,$ and so on. 
	\medskip
{\bf Remark.} {\it If $M$ is an indecomposable $3$-dimensional $\Lambda$-module, then its annihilator is equal to $U(a,b,c)$ for some $(a,b,c)\neq 0$ and 
$M$ considered as a
$\Lambda/U(a,b,c)$-module is either the unique indecomposable projective
$\Lambda/U(a,b,c)$-module} (and then a local module, thus isomorphic to
$M(a,b,c)$) {\it or the unique indecomposable injective
$\Lambda/U(a,b,c)$-module} (and then a module with simple socle, thus isomorphic to
$D(M(a,b,c))$). 

   \medskip
Proof of the Proposition and the Remark. First, let us assume that $M$ is local.
According to (2.6) and (1.4), we know that $M \simeq M(a\col b\col c)$ for some
$(a\col b\col c) \in \Bbb P^2$ and that these modules are pairwise non-isomorphic.
As representatives of the elements of $\Bbb P^2$, we choose (as usual) the triples
$(c_1,c_2,c_3)$ with $c_i = 1$ for some $i$ and $c_j = 0$ for $j<i$. Clearly, there
are the seven cases (1) to (7) as listed above. 
It remains to choose in every case a basis 
$\Cal B(a,b,c) = \{v,v_1,v_2\}$ of $M(a,b,c)$. Recall that
$M(a,b,c) = \overline\Lambda/(a\col b\col c)$ is a factor module of $\Lambda$
and $\overline\Lambda$ has the basis $\{1,x,y,z\}$. 
We choose as elements of $\Cal B(a,b,c)$ the residue class $v = \overline 1$ as well as 
two of the three residue classes
$\overline x, \overline y, \overline z$, namely 
$v_1 = \overline x$ if $a = 0$ and $v_1 = \overline y$ otherwise, and then
$v_2 = \overline y$ in case $(a,b,c) = (0,0,1)$ and $v_2 = \overline z$ otherwise. 
(We should remark that the vertices and the arrows of the diagram are those 
of the coefficient quiver $\Gamma(M(a,b,c),\Cal B(a,b,c))$ as considered in [R],
and the solid arrows focus the attention to a spanning tree.)

Second, assume that $M$ is not local. Since $M$ is an indecomposable module of
length 3 and Loewy length 2, it follows that $M$ has simple socle, thus $D(M)$
is local and therefore of the form (1) to (7).

Finally, $M$ and $D(M)$ have the same annihilator, this is a 3-dimensional ideal,
thus of the form $U(a,b,c)$. The 3-dimensional local algebra $\Lambda/U(a,b,c)$
has a unique 3-dimensional local module, this is the indecomposable projective
$\Lambda/U(a,b,c)$-module, and dually, it has a unique 3-dimensional module
with simple socle, this is the unique  indecomposable injective
$\Lambda/U(a,b,c)$-module.
This completes the proof.
$\s$

	\bigskip\medskip
{\bf Reference.}
	\medskip
\item{[AIS]} L. L. Avramov, S. B. Iyengar, L. M. \c Sega.
  Free resolutions over short local rings. J. London Math. Soc. (2) 78 (2008),
  459--476.
\item{[CV]} L. W. Christensen, O. Veliche.
 Acyclicity over local rings with radical cube zero. Illinois J. Mathematics. 51
 (2007), 1439--1454.
\item{[J\c S]} D. A. Jorgensen, L. M. \c Sega. Independence of the total reflexivity
 conditions for modules. Algebras and Representation Theory 9 (2006), 217--226.
\item{[R]} C\. M\. Ringel. Exceptional modules are tree modules.
    Lin. Alg. Appl. 275--276 (1998). 471--493.
\item{[RZ1]} C. M. Ringel, P. Zhang: Gorenstein-projective and
 semi-Gorenstein-projective modules. To appear.
 arXiv:1808.01809v3.
\item{[RZ2]} C. M. Ringel, P. Zhang. Gorenstein-projective modules over short
 local algebras. In preparation.

\bigskip

\baselineskip=1pt
{\rmk
C. M. Ringel\par
Fakult\"at f\"ur Mathematik, Universit\"at Bielefeld \par
POBox 100131, D-33501 Bielefeld, Germany  \par
ringel\@math.uni-bielefeld.de
\medskip

P. Zhang \par
School of Mathematical Sciences, Shanghai Jiao Tong University \par
Shanghai 200240, P. R. China.\par
pzhang\@sjtu.edu.cn}

\bye